\documentclass[11pt]{article}
\usepackage{amsfonts,latexsym,amsmath,amscd,geometry}
\usepackage{amssymb}
\usepackage{latexsym}

\newcommand \nc{\newcommand}
\newtheorem{theorem}{Theorem}[section]
\newtheorem{lemma}[theorem]{Lemma}

\newtheorem{corollary}[theorem]{Corollary}

\newtheorem{remark}[theorem]{Remark}

\nc{\ba}{\begin{array}}\nc{\ea}{\end{array}}
\nc{\be}{\begin{eqnarray}}\nc{\ee}{\end{eqnarray}}
\nc{\beq}{\begin{equation}}\nc{\eeq}{\end{equation}}
\nc{\bex}{\begin{eqnarray*}}\nc{\eex}{\end{eqnarray*}}
\nc{\btm}{\begin{theorem}} \nc{\etm}{\end{theorem}}
\nc{\blm}{\begin{lemma}} \nc{\elm}{\end{lemma}}
\nc{\R}{\mathbb{R}} \nc{\va}{\varepsilon} \nc{\ls}{\limits}

\def\pf{\noindent{\bf Proof.\quad}}

\newcommand \qed {\hfill $\Box$}

\begin{document}

\title{Remarks on approximate harmonic maps in dimension two}
\author{Changyou Wang\footnote{Department of Mathematics, Purdue University, West Lafayette, IN 47907, USA}}
\maketitle

\begin{abstract}
For the class of approximate harmonic maps $u\in W^{1,2}(\Sigma,N)$ from a closed Riemmanian surface $(\Sigma,g)$ to a compact Riemannian
manifold $(N, h)$, we show that (i) the so-called energy identity holds for weakly convergent
approximate harmonic maps $\{u_n\}:\Sigma\to N$, with tension fields $\tau(u_n)$
bounded in the Morrey space $M^{1,\delta}(\Sigma)$ for some
$0\le\delta<2$;  and (ii) if an approximate harmonic map $u$ has tension field 
$\tau(u)\in L\log L(\Sigma)\cap M^{1,\delta}(\Sigma)$ for some
$0\le\delta<2$, then $u\in W^{2,1}(\Sigma, N)$. Based on these estimates, we further establish
the bubble tree convergence, referring to
energy identity both $L^{2,1}$ of gradients and $L^1$-norm of hessians and the 
oscillation convergence, for a weakly convergent sequence of
approximate harmonic maps $\{u_n\}$, with tension fields $\tau(u_n)$ uniformly bounded in
$M^{1,\delta}(\Sigma)$ for some $0\le\delta<2$ and uniformly integrable 
in $L\log L(\Sigma)$. 
\end{abstract}

\section{Introduction}
\setcounter{equation}{0}
\setcounter{theorem}{0}

Harmonic maps form one of the most important classes of geometric partial differential equations.
They are critical points of the Dirichlet energy functional
$$E(u)=\frac12\int_\Sigma |\nabla u|^2\,dv_g$$
for maps $u\in W^{1,2}(\Sigma, N)$ between two Riemannian manifolds $(\Sigma^m, g)$ and $(N^l, h)$, which is isometrically
embedded into some euclidean space $\R^L$. The equation for harmonic maps is
\begin{equation}\label{HM1}
\Delta u+A(u)(\nabla u,\nabla u)=0,
\end{equation}
where $A(\cdot)(\cdot,\cdot)$ is the second fundamental form of $N$ in $\R^L$. 

Throughout this paper, we assume that {\it $m={\rm{dim}}
(\Sigma)=2$, i.e. $(\Sigma,g)$ is a Riemannian surface.}  
Since the Dirichlet energy is invariant under conformal transformations of $\Sigma$, 
solutions to the harmonic map equation are subject to concentration compactness phenomena and hence
are of particular interest.

For a positive integer $k$ and $1\le p\le\infty$, recall that the Sobolev space $W^{k,p}(\Sigma,N)$ is defined by
$$W^{k,p}(\Sigma, N)=\Big\{u\in W^{k,p}(\Sigma, \R^L): \ u(x)\in N \ {\rm{a.e.}}\ x\in \Sigma\Big\}.$$
In a seminal work \cite{Helein}, H\'elein established the smoothness of any solution $u\in W^{1,2}(\Sigma, N)$ to \eqref{HM1}. Subsequently, Evans \cite{evans} and Bethuel \cite{bethuel} have obtained the partial regularity theorem on stationary harmonic maps in higher dimensions $m\ge 3$, see also \cite{CWY}
and \cite{RS} for alternate proofs. See also \cite{Lin-Wang2} for more references.
The original ideas
by \cite{Helein} are based on rewriting \eqref{HM1} into a form in which the nonlinearity enjoys
finer algebraic structures under  Coulomb gauge frames on $N$ so that some compensated regularity
can be achieved by employing harmonic analysis techniques from Hardy or Lorentz spaces. The ideas
of \cite{Helein} were significantly extended by Rivi\`ere \cite{riviere}, in which the continuity was shown
for any weak solution $u\in W^{1, 2}(B_1^2,\R^L)$ to the equation:
\begin{equation}\label{HM2}
-\Delta u=\Omega\cdot \nabla u, 
\end{equation}
for any given $\Omega\in L^2(B_1^2, so(L)\otimes\R^2)$. We point out that \eqref{HM2} is more general than \eqref{HM1}, which includes harmonic map equations,
the $H$-surface equation, and critical points of any conformally invariant elliptic Lagrangian that is quadratic in the gradient.

Because of the conformal invariance of the Dirichlet energy in dimension two, weakly convergent sequences of harmonic maps $\{u_n\}\subset W^{1,2}(\Sigma, N)$ may not converge strongly in 
$W^{1,2}(\Sigma, N)$. This is the so-called bubbling phenomena, i.e., the energy can concentrate 
at a set of finitely many points, at where  a nontrivial harmonic map from $\mathbb S^2$ to $N$, called 
a {\it bubble}, can generate. It has been an important question to establish the energy identity that
accounts the total energy loss by that of a number of bubbles during the bubbling process. Such an
energy identity has first been proven by Jost \cite{jost} and Parker \cite{Parker} for any sequence
of harmonic maps. 

Motivated by the studies, such as the long time dynamics and behavior near singularities, of the heat flow of harmonic maps from surfaces initiated by Struwe \cite{struwe},
people have also been interested in the bubbling phenomena for approximate harmonic maps
from surfaces. Recall that a map $u\in W^{1,2}(\Sigma, N)$ is called 
an approximate harmonic map  with tension field $\tau(=\tau(u)):\Sigma\to T_u N$, if
\begin{equation}\label{HM3}
\Delta u+A(u)(\nabla u, \nabla u)=\tau.
\end{equation}
In particular, an approximate harmonic map $u$ becomes a harmonic map if its tension field $\tau(u)=0$.
While for any function $\tau\in L^1(\Sigma,\R^L)$, we can show that any critical point $u\in W^{1,2}(\Sigma,N)$ to 
$$E(u,\tau):=\int_\Sigma \big(\frac12|\nabla u|^2+\langle\tau, u\rangle \big)\,dv_g$$
is an approximate harmonic map with tension field $\tau(u)=P_u(\tau)$, where $P_y:\R^L\to T_y N$,
$y\in N$, is an orthogonal projection. The question concerning approximate harmonic maps seeks
general sufficient conditions on their tension fields to guarantee the energy identity or the so-called
bubble tree convergence, i.e., the oscillation convergence.  There have been many works in this direction.
When $\tau(u_n)$ is bounded in $L^2(\Sigma)$, the energy identity has been proved
by Qing \cite{Qing} for $N=\mathbb S^{L-1}$ and Ding-Tian \cite{Ding-Tian} and Wang \cite{wang} for
general target manifolds $N$; and the bubble tree convergence has been established by
Qing-Tian \cite{Qing-Tian} and Lin-Wang \cite{Lin-Wang}. We remark that the space
$L^2(\Sigma)$ for $\tau(u)$ is not conformally invariant. While $L^1(\Sigma)$ for $\tau(u)$ is conformally invariant, an example by Parker \cite{Parker} showed that the energy identity fails if $\tau(u_n)$
is merely assumed to be bounded in $L^1(\Sigma)$.  Therefore, an interpolation space between
$L^1(\Sigma)$ and $L^2(\Sigma)$ for $\tau(u)$ seems necessary for the energy identity. In this
direction, the energy identity has been proved by Lin-Wang \cite{Lin-Wang1} when $N=\mathbb S^{L-1}$
and $\tau(u_n)$ is bounded in $L^p(\Sigma)$ for some $1<p\le 2$, by Li-Zhu \cite{Li-Zhu} for general
target manifolds $N$ when $\tau(u_n)$ is bounded in $L^p(\Sigma)$ for some $p\ge\frac65$ and
by Luo \cite{Luo} under the condition $\big\|\tau(u_n)\big\|_{L^2(B_r(x)\setminus B_{\frac{r}2}(x))}\lesssim r^{-a}$ for some $a\in (0,1)$ and any $r\in (0,1)$ and $x\in\Sigma$. Most recently, 
in a very interesting article \cite{WWZ} Wang-Wei-Zhang have obtained  that the energy
identity and bubble tree convergence holds for approximate harmonic maps $u_n$ whenever $\tau(u_n)$
is bounded in $L^p(\Sigma)$ for some $1<p\le 2$ by developing some delicate estimates based on H\'elein's
original approach in \cite{Helein}. We would also like to mention that there have been many
interesting works on the bubbling analysis of solutions to the equation \eqref{HM2} by Laurain-Rivi\`ere
\cite{Lau-Riv}, and by Sharp-Topping \cite{Sharp-Topping} and Lamm-Sharp \cite{Lamm-Sharp} on 
perturbed version of \eqref{HM2}: 
\begin{equation}\label{HM4}
-\Delta u=\Omega\cdot\nabla u+f,
\end{equation}
for some non-homogeneous term $f$. 
Among other results, it was shown by \cite{Lamm-Sharp} that (i) in dimension $2$, a global $W^{2,1}$-estimate holds for an approximate harmonic map $u$ when $\tau(u)\in L\log L(\Sigma)$ and $|\mathcal H(u)|^\frac12\in L^{2,1}(\Sigma)$, here $\mathcal H(u)$ is the Hopf differential of $u$ (see \cite{Helein} ($m=2$) and Lin-Rivi\`ere \cite{Lin-Riviere} ($m\ge 3$) for harmonic maps into $N=\mathbb S^{L-1}$); and (ii) an energy identity
holds when $\tau(u_n)$ is bounded in $L\log L(\Sigma)$ and $|\mathcal H(u_n)|^\frac12(\Sigma)$ is bounded in $L^{2,1}(\Sigma)$,
and the bubble tree convergence holds if $\tau(u_n)$ is bounded in $L^p(\Sigma)$ for some $p>1$ and
$|\mathcal H(u_n)|^\frac12$ is uniformly integrable in $L^{2,1}(\Sigma)$.

In this note, we aim to improve the results on the global $W^{2,1}$-estimate of approximate harmonic maps
by \cite{Lamm-Sharp}, and the energy identity and bubble tree convergence of approximate harmonic maps
by \cite{WWZ} by relaxing the conditions on the tension fields.

Before stating our results, we need to introduce several definitions.
For a bounded domain $U\subset\R^2$, $1\le p<+\infty$, and $0\le\delta\le 2$, we define the Morrey ($p,\delta$)-space on $U$ by
$$M^{p,\delta}(U):=\Big\{f\in L^p_{\rm{loc}}(U):
\big\|f\big\|^p_{M^{p,\delta}(U)}\equiv\sup_{B_r(x)\subset U} r^{\delta-2}\int_{B_r(x)}|f|^p<\infty\Big\}.$$
It is readily seen that 
$$L^\infty(U)=M^{p,0}(U)\subset M^{p,\delta_1}(U)\subset M^{p,\delta_2}\subset L^p(U)
=M^{p,2}(U)
$$ 
for any $0\le\delta_1\le\delta_2\le 2$. 

The space $L\log L(U)$ is defined by
$$L\log L(U)=\Big\{f\in L^1(U): \big\|f\big\|_{L\log L(U)}=\int_{U} |f|\log(2+|f|)<+\infty\Big\}.$$
For any $1<p<\infty$ and $1\le q\le \infty$, we denote by $L^{p,q}(U)$
as the Lorentz ($p,q$)-space, see Ziemer \cite{Ziemer} for its definition. 
In this paper, we will mainly work
with the spaces $L^{2,1}(U)$ and $L^{2,\infty}(U)$, which can be defined as
follows (cf. \cite{Lau-Riv}):  
$$L^{2,1}(U)=\Big\{f\in L^1_{\rm{loc}}(U)\ |\ \|f\|_{L^{2,1}(U)}\equiv\int_0^\infty \lambda_f^\frac12(t)\,dt<\infty\Big\},
$$
$$L^{2,\infty}(U)=\Big\{f\in L^1_{\rm{loc}}(U)\ |\ \|f\|_{L^{2,\infty}(U)}\equiv\sup_{t>0} 
\ t\lambda_f^\frac12(t)\,dt<\infty\Big\},
$$
where $\lambda_f(t)=\big|\{x\in U: |f(x)|\ge t\}\big|$.  It is well-known (cf. \cite{Helein}) that
$(L^{2,1}(U))^*=L^{2,\infty}(U)$.

Our first result concerns the global $W^{2,1}$-estimate on approximate harmonic maps in dimension two,
which is stated as follows.

\begin{theorem}\label{21-estimate1}
If the tension field $\tau(u)$ of an approximate harmonic map $u\in W^{1,2}(U, N)$ satisfies 
$\tau(u)\in L\log L(U)\cap M^{1,\delta}(U)$ for some $0\le\delta<2$, then $u\in W^{2,1}_{\rm{loc}}(U,N)$.
Moreover, for any compact subset $K\subset U$, it holds that
\begin{equation}\label{21-estimate2}
\big\|\nabla u\big\|_{L^{2,1}(K)}+\big\|\nabla^2 u\big\|_{L^{1}(K)}\le C\big(\delta, K, \|\nabla u\|_{L^2(U)}, \|\tau(u)\|_{L\log L(U)},
\|\tau(u)\|_{M^{1,\delta}(U)}\big).
\end{equation}
\end{theorem}

Recall that if $f\in L^p(U)$ for some $1<p\le 2$, then $f\in L\log L(U)\cap M^{1,\frac{2}{p}}(U)$, and
$$\big\|f\|_{L\log L(U)}+\big\|f\big\|_{M^{1,\frac{2}{p}}(U)}\le C(p,U)\big\|f\big\|_{L^p(U)}.$$
Thus the conclusion of Theorem \ref{21-estimate1} holds whenever $\tau(u)\in L^p(U)$ for some
$1<p\le 2$.

Our second result concerns both the energy identity and bubble tree convergence of approximate
harmonic maps in dimension two. 
\begin{theorem}\label{bubbling1}
Let $\{u_n\}\subset W^{1,2}(\Sigma,N)$ be a sequence of approximate harmonic maps, with tension fields
$\tau(u_n)$,  converging weakly
to a map $u\in W^{1,2}(\Sigma,N)$ such that
\begin{equation}\label{tension-bound}
\sup_{n}\big\|\tau(u_n)\big\|_{M^{1,\delta}(\Sigma)}<\infty
\end{equation}
for some $0\le\delta<2$. Let $\tau\in L^1(\Sigma)$ be such that
$\tau(u_n)\rightharpoonup \tau$ in $L^1(\Sigma)$. Then\\
(i) $u$ is an approximate harmonic map, with tension field $\tau\in  M^{1,\delta}(\Sigma)$, and
$u\in C^{2-\delta}(\Sigma)\cap W^{1,q}(\Sigma)$ for any $2\le q<q(\delta)$. Here
$q(\delta)=\frac{\delta}{\delta-1}$ if $1<q<2$; and $q(\delta)=\infty$ if $0\le \delta\le 1$.
\\
(ii) there exist  $m$-nontrivial harmonic maps $\{\omega_i\}_{i=1}^m\subset C^\infty(\mathbb S^2,N)$, $\{x_n^i\}_{i=1}^m\subset \Sigma,\ \{r_n^i\}_{i=1}^m\subset\R_+$, with 
$\displaystyle\lim_{n\rightarrow\infty}r_n^i=0$ ($i=1,\cdots, m$), and
\begin{equation}\label{scale-separation}
\lim_{n\rightarrow\infty}\Big\{\frac{r_n^i}{r_n^j}, \frac{r_n^j}{r_n^i}, \frac{|x_n^i-x_n^j|}{r_n^i+r_n^j}\Big\}=\infty,
\ 1\le i<j\le m, 
\end{equation}
such that
\begin{equation}\label{bubble2.0}
\Big\|u_n(\cdot)-u(\cdot)-\sum_{i=1}^m \big(\omega_i(\frac{\cdot-x_n^i}{r_n^i})-\omega_i(\infty)\big)\Big\|
_{W^{1,2}(\Sigma)}=0.
\end{equation}
In particular, 
\begin{equation}\label{bubbling2}
\lim_{n\rightarrow\infty}\int_\Sigma|\nabla u_n|^2
=\int_\Sigma|\nabla u|^2+\sum_{i=1}^m \int_{\mathbb S^2}|\nabla\omega_i|^2.
\end{equation}
(iii) if, in addition, $\tau(u_n)$ is  uniformly integrable in $L\log L$, i.e.,
\begin{equation}\label{equi-cont}
\lim_{|E|\rightarrow 0}\sup_{n}\big\|\tau(u_n)\big\|_{L\log L(E)}=0,
\end{equation}
then 
\begin{eqnarray}\label{bubbling3}
\lim_{n\rightarrow\infty}\big\|\nabla u_n\big\|_{L^{2,1}(\Sigma)}
=\big\|\nabla u\big\|_{L^{2,1}(\Sigma)}+\sum_{i=1}^m \big\|\nabla\omega_i\big\|_{L^{2,1}(\mathbb S^2)}, 
\end{eqnarray}
and
\begin{equation}\label{oscillation-convergence}
\Big\|u_n(\cdot)-u(\cdot)-\sum_{i=1}^m \big(\omega_i(\frac{\cdot-x_n^i}{r_n^i})-\omega_i(\infty)\big)\Big\|_{L^\infty(\Sigma)}=0.
\end{equation}
\end{theorem}

As remarked above, if $f_n$ is bounded in $L^p(\Sigma)$ for some $1<p\le 2$, then
$f_n$ is bounded in $M^{1,\frac{2}{p}}(\Sigma)$ and is uniformly integrable in $L\log L$.
Therefore, the conclusions of Theorem \ref{bubbling1} hold under the condition that $\tau(u_n)$
remains to be bounded in $L^p(\Sigma)$ for some $1<p\le 2$.

During the proof of Theorem 1.2 (iii) in section 3 below, we actually show that there is no $W^{2,1}$-norm
concentration of $u_n$ in the neck region.  As a byproduct of this conclusion, we can prove the following
Corollary. 

\begin{corollary}\label{W21-convergence}
Assume $\{u_n\}\subset W^{1,2}(\Sigma,N)$ is a sequence of approximate harmonic maps,
with tension fields $\tau(u_n)$,  converging weakly
to an approximate harmonic map $u\in W^{1,2}(\Sigma,N)$, with tension field $\tau$, such that
\begin{equation}\label{tension-bound1}
\sup_{n}\Big(\big\|\tau(u_n)\big\|_{L\log L(\Sigma)}+\big\|\tau(u_n)\big\|_{M^{1,\delta}(\Sigma)}\Big)<\infty
\end{equation}
for some $0\le\delta<2$. If, in addition, 
\begin{equation}\label{LlogL}
\lim_{n\rightarrow\infty}\big\|\tau(u_n)-\tau\big\|_{L\log L(\Sigma)}=0,
\end{equation}
there exist finitely many nontrivial harmonic maps $\{\omega_i\}_{i=1}^m\subset C^\infty(\mathbb S^2,N)$ such that
\begin{equation}\label{bubbling4}
\lim_{n\rightarrow\infty}\big\|\nabla^2 u_n\big\|_{L^{1}(\Sigma)}
=\big\|\nabla^2 u\big\|_{L^{1}(\Sigma)}+\sum_{i=1}^m \big\|\nabla^2\omega_i\big\|_{L^{1}(\mathbb S^2)}.
\end{equation}
\end{corollary}

We would like to remark that \eqref{bubbling4} was previously proven by Lamm-Sharp \cite{Lamm-Sharp}
for approximate harmonic maps under a stronger condition that the Hopf differential $|\mathcal H(u_n)|^\frac12$
is uniformly integrable in $L^{2,1}(\Sigma)$, and $\tau(u_n)$ is bounded  in $W^{1,2}(\Sigma)$.

The argument to prove both theorems is built upon the interesting approach outlined by Wang-Wei-Zhang \cite{WWZ} together with some new observations.

The paper is organized as follows. In section 2, we will provide all the necessary estimates, in particular,
the holomorphic approximation of Hope differential. In section 3, we will prove both Theorem 1.1, Theorem
1.2, and Corollary 1.3.
 
\section{Apriori estimates on approximate harmonic maps}
\setcounter{equation}{0}
\setcounter{theorem}{1}
This section is devoted to several apriori estimates on approximate harmonic maps under various
conditions on their tension fields. It includes the crucial estimate of a holomorphic approximation
of the Hopf differential of an approximate harmonic map $u$ with tension field in 
$\tau(u)\in M^{1,\delta}(\Sigma)$,  which is an improvement of Proposition 4.1 of Wang-Wei-Zhang \cite{WWZ} 
where $\tau(u)\in L^p(\Sigma)$ for some $p>1$ is assumed.

\subsection{$\bar{\partial}$-equation for approximate harmonic maps}
Let $B_r(x)\subset\R^2$ denote the ball with center at $x\in\R^2$ and radius $r>0$, and $B_r=B_r(0)$.
Assume that $U\subset B_2$ is an open set and $u:U\to N$ is an approximate harmonic map.
From H\'elein's argument of enlarging the target manifolds, we may assume that
$(N,h)\subset\R^L$ is parallelized and there exists a global orthonormal frame
$\{e_\alpha\}_{\alpha=1}^l$, $l={\rm{dim}}(N)$, of the pull-back tangent bundle
$u^*TN$ on $U\subset\R^2$ such that
\begin{equation}\label{col-frame}
\begin{cases}
{\rm{div}}\langle\nabla e_\alpha, e_\beta\rangle =0 \ {\rm{in}}\ U, \ 1\le\alpha,\beta\le l,\\
\displaystyle\sum_{\alpha=1}^l\int_{U}|\nabla e_\alpha|^2\le C\int_U|\nabla u|^2.
\end{cases}
\end{equation}
Moreover, it follows from H\'elein \cite{Helein} that $\nabla e_\alpha\in L^{2,1}(U)$ and
\begin{equation}\label{frame-est}
\sum_{\alpha=1}^l\big\|\nabla e_\alpha\big\|_{L^{2,1}(U)}\le C\int_{U}|\nabla u|^2.
\end{equation}
Denote 
$$\frac{\partial}{\partial z}=\frac{\partial}{\partial x}-\frac{\partial}{\partial y} i,\
\frac{\partial}{\partial \bar{z}}=\frac{\partial}{\partial x}+\frac{\partial}{\partial y} i.
$$
Then any approximate harmonic map $u\in W^{1,2}(U,N)$ can be written as
\begin{equation}\label{approx-hm1}
\frac{\partial}{\partial {\bar z}}\langle\frac{\partial u}{\partial z}, e_\alpha\rangle
=\langle \frac{\partial u}{\partial z}, e_\beta\rangle \langle \frac{\partial e_\alpha}{\partial {\bar z}}, e_\beta\rangle
+\langle\tau(u), e_\alpha\rangle\ \ {\rm{in}}\ \ U.
\end{equation}
Set $$\displaystyle G(u)=\big(\langle\frac{\partial u}{\partial z}, e_\alpha\rangle\big)_{1\le\alpha\le l},
\ T(u)=\big(\langle \tau(u), e_\alpha\rangle\big)_{1\le\alpha\le l}: U\to\mathbb C^l,$$
and the connection matrix on $u^*TN$ 
$$\omega(u)=\big(\langle\frac{\partial e_\alpha}{\partial \bar z}, e_\beta\rangle\big)_{1\le\alpha,\beta\le l}:
U\to\mathbb C^{l\times l}.$$
Then \eqref{approx-hm1} can be rewritten as
\begin{equation}\label{approx-hm2}
\frac{\partial}{\partial\bar z} (G(u))=\omega(u) G(u)+T(u) \ \ {\rm{in}}\ \ U.
\end{equation}

Define the linear operator $\mathcal T: L^\infty\big(\mathbb C, \mathbb C^{l\times l}\big)
\to L^\infty\big(\mathbb C, \mathbb C^{l\times l}\big)$ by
\begin{equation}\label{t-operator}
(\mathcal{T}A)(z)=\int_{\mathbb C}\frac{\omega(u(\xi))\chi_{U}(\xi) A(\xi)}{\pi(z-\xi)}\,d\xi, \ z\in \mathbb C,
\end{equation}
where $\chi_U$ is the characteristic function of $U$.
It is readily seen that
\begin{equation}\label{d-bar-eqn}
\frac{\partial}{\partial\bar z} A=\omega(u) A \ \ {\rm{in}}\ \  U.
\end{equation}
Since $\displaystyle\frac{1}{z}\in L^{2,\infty}(\mathbb C)$ has
$$\big\|\frac{1}{z}\big\|_{L^{2,\infty}(\mathbb C)}\le \sqrt{\pi},$$
we have that for any $A\in L^\infty\big(\mathbb C, \mathbb C^{l\times l}\big)$, it holds that
\begin{equation}\label{t-operator-est}
\big\|\mathcal{T} A\big\|_{L^\infty(\mathbb C)}\le \frac{1}{\pi}\big\|\omega(u)\big\|_{L^{2,1}(U)}
\big\|\frac{1}{z}\big\|_{L^{2,\infty}(\mathbb C)}\big\|A\big\|_{L^\infty(\mathbb C)}\le C\epsilon_0 \big\|A\big\|_{L^\infty(\mathbb C)}.
\end{equation}
Based on \eqref{t-operator-est}, H\'elein \cite{Helein} has shown  
\begin{lemma}\label{conjugate} 
There exists $\epsilon_0>0$ such that if $\displaystyle\|\nabla u\big\|_{L^{2}(U)}
\le\epsilon_0$, then there exists a unique solution $B\in L^\infty\big(\mathbb C, \mathbb C^{l\times l}\big)$ of
\begin{equation}\label{t-equation}
B-\mathcal{T} B=\mathbb I_l\ \ {\rm{in}}\ \ \mathbb C, 
\end{equation}
such that
\begin{equation}\label{invertible}
\big\|B-\mathbb I_l\big\|_{L^\infty(\mathbb C)}\le \frac12, \ {\rm{and}}\  B^T B=\mathbb I_l.
\end{equation}
Here $\mathbb I_l$ is the identity matrix of order $l$.
\end{lemma}
It follows from \eqref{invertible} that $B$ is invertible and solves 
\begin{equation}\label{harmonic1}
\frac{\partial}{\partial\bar z} B=\omega(u) B\ \ {\rm{in}}\ \ U.
\end{equation}
In particular, $B^{T}$ solves
\begin{equation}\label{harmonic2}
\frac{\partial}{\partial\bar z} B^{T}=-\omega(u) B^{T}\ \ {\rm{in}}\ \ U.
\end{equation}
Combining \eqref{approx-hm2} with \eqref{harmonic2}, we obtain that
\begin{equation}\label{approx-hm3}
\frac{\partial}{\partial\bar z}(B^{T}G(u))=B^{T} T(u) \ \ {\rm{in}}\ \  U.
\end{equation}
\subsection{H\"older continuity estimate of approximate harmonic maps}
Employing  \eqref{approx-hm3}, we can show the H\"older continuity for approximate
harmonic maps.
Let $\epsilon_0>0$ be the constant given by Lemma \ref{conjugate} 
so that \eqref{approx-hm3} holds.
Set $G_1, G_2:U\to\mathbb C^{l}$ by
\begin{equation}\label{decom1}
\begin{cases}
G_1(z)=\int_{U}\frac{B^{T}(\xi)T(u)(\xi)}{\pi(z-\xi)}\,d\xi,\ z\in U,\\
G_2(z)=B^{T}G(u)-G_1(z), \  z\in U.
\end{cases}
\end{equation}
Thus it holds 
\begin{equation}\label{decom2}
B^{T}G(u)=G_1+G_2\ {\rm{in}}\ U.
\end{equation}
We can prove the following Lemma.

\begin{lemma}\label{epsilon-regularity1} There exists $\epsilon_0>0$ such that
if $u\in W^{1,2}(U, N)$ is an approximate harmonic map, satisfying $\displaystyle\int_{U} |\nabla u|^2\le\epsilon_0^2$
and $\tau(u)\in M^{1,\delta}(U)$ for some $1<\delta<2$, 
then \\
(i) $G_1\in M^{\frac{\delta}{\delta-1},\delta}_*(U)$, and
\begin{equation}\label{G1-morrey-est}
\big\|G_1\big\|_{M^{\frac{\delta}{\delta-1},\delta}_*(U)}
\le C\big\|\tau(u)\big\|_{M^{1,\delta}(U)}.
\end{equation}
Here
$$\big\|f\big\|_{M^{p,\delta}_*(U)}^{p}:=\sup_{B_r(x)\subset U}
\Big\{ r^{\delta-2}\sup_{\lambda>0} \lambda^{p}\big|\{y\in B_r(x): \ |f(y)|\ge\lambda\}\big|\Big\},
\ 1\le p<\infty.
$$
(ii) $u\in C^{2-\delta}(U)$, and
\begin{equation}\label{holder-estimate}
\big[u\big]_{C^{2-\delta}(B_{\frac{R}2}(x_0))}\le C(\delta, R) \big(\epsilon_0+\big\|\tau(u)\big\|_{M^{1,\delta}(U)}\big)
\end{equation}
for any ball $B_R(x_0)\subset U$.
\end{lemma}
\pf  It follows from \eqref{approx-hm3} that $G_2$ is holomorphic in
$U$, i.e., $\frac{\partial G_2}{\partial \bar z}=0$ in $U$. Since $B^{T}$ satisfies
$\frac23\le\displaystyle\big\|B^{T}\big\|_{L^\infty(U)}\le 2$, it follows from \eqref{decom1} and \eqref{decom2}
that
\begin{equation}\label{reisz1}
|\nabla u(z)|\le C\Big(|G_2(z)|+\int_{U}\frac{|\tau(u)(\xi)|}{|z-\xi|}\,d\xi\Big)
=C\big(|G_2(z)|+\mathcal {I}_1(\tau(u))(z)\big),\ z\in U, 
\end{equation}
where $\mathcal I_1$ is the Riesz potential of order 1 on $U$:
$$\mathcal I_1(f)(z)=\int_{U}\frac{|f(\xi)|}{|z-\xi|}\,d\xi, \  z\in U, \ f\in L^1(U).$$
Since $\tau(u)\in M^{1,\delta}(U)$, applying Adams' theorem on the estimate of Riesz potential in Morrey spaces (cf. \cite{adams}) we 
have that $\mathcal I_1(\tau(u))\in M^{\frac{\delta}{\delta-1}, \delta}_*(U)$, and
\begin{equation}\label{morrey-est}
\big\|\mathcal I_1(\tau(u))\big\|_{M^{\frac{\delta}{\delta-1}, \delta}_*(U)}
\le C\big\|\tau(u)\big\|_{M^{1,\delta}(U)}.
\end{equation}
This yields (\ref{G1-morrey-est}).

It follows from (\ref{reisz1}) and (\ref{morrey-est})
that for any given ball $B_{R}(x_0)\subset U$, if $B_r(x)\subset B_{\frac{R}2}(x_0)$
then it holds
$$r^{\delta-2}\big\|\nabla u\big\|_{L^{\frac{\delta}{\delta-1}}_*(B_r(x))}^{\frac{\delta}{\delta-1}}
\le C\big(\|G_2\|_{L^\infty(B_{\frac{R}2}(x_0))}^{\frac{\delta}{\delta-1}}r^\delta+\big\|\tau(u)\big\|_{M^{1,\delta}(U)}\big).
$$
Since $G_2$ is holomorphic, by the mean value property we have
\begin{eqnarray*}
\|G_2\|_{L^\infty(B_{\frac{R}2}(x_0))}&\le& CR^{-2}\int_{B_R(x_0)}|G_2|
\le CR^{-2}\int_{U}(|\nabla u|+|G_1|)\\
&\le& CR^{-2}\int_{U}(|\nabla u|+|\tau(u)|),
\end{eqnarray*}
where we have used Young's inequality to estimate $G_1$: 
$$\|G_1\|_{L^1(U)}\le \|\mathcal I_1(\tau(u))\|_{L^1(U)}\le C\|\tau(u)\|_{L^1(U)}.$$
Therefore we obtain
\begin{eqnarray*}
r^{\delta-2}\big\|\nabla u\big\|_{L^{\frac{\delta}{\delta-1}}_*(B_r(x))}^{\frac{\delta}{\delta-1}}
\le C(R)\big(\|\nabla u\|_{L^2(U)}+ \|\tau(u)\|_{M^{1,\delta}(U)}\big),
\end{eqnarray*}
so that
$$
r^{\frac{\delta}{\delta-1}-2}\big\|\nabla u\big\|_{L^{\frac{\delta}{\delta-1}}_*(B_r(x))}^{\frac{\delta}{\delta-1}}
\le C(R)\big(\|\nabla u\|_{L^2(U)}+\|\tau(u)\|_{M^{1,\delta}(U)}\big) r^{\frac{\delta}{\delta-1}(2-\delta)},
\ \forall\ B_r(x)\subset B_{\frac{R}2}(x_0).
$$
This, with Morrey's decay Lemma \cite{morrey}, 
implies that $u\in C^{2-\delta}(B_{\frac{R}2}(x_0))$ and
\eqref{holder-estimate} holds.
\qed

\begin{remark}\label{2-morrey-estimate} {\rm For $1<\delta<2$, since $L^2(B_r(x))\subset 
L^{\frac{\delta}{\delta-1}}_*(B_r(x))$, it is not hard to see that \eqref{G1-morrey-est} implies
that}
\begin{equation}\label{G1-morrey-est1}
\big\|G_1\big\|_{M^{2,2\delta-2}(U)}\le C\big\|G_1\big\|_{M^{\frac{\delta}{\delta-1},\delta}_*(U)}
\le C\big\|\tau(u)\big\|_{M^{1,\delta}(U)}.
\end{equation}
\end{remark}
\subsection{Apriori estimates of $L^2$ and $L^{2,1}$-norms of $\nabla u$ in the neck region}

We now want to apply Lemma 2.3 to estimate both $L^2$ and $L^{2,1}$-norms of $\nabla u$ in the neck region.
For this purpose, let $U=B_1\setminus B_r$, $0<r\le \frac14$, and
$u: B_1\setminus B_r\to N$ be an approximate harmonic map satisfying
the conditions of Lemma 2.3. Using the same notations as in the proof of Lemma 2.3,
we can write the Laurent series of $G_2$ on $B_1\setminus B_r$ as
\begin{equation}\label{g2-form}
G_2(z)=\sum_{n=-\infty}^\infty a_n z^n, \ z\in B_1\setminus B_r;
\end{equation}
while
\begin{equation}\label{g1-form}
G_1(z)=\int_{B_1\setminus B_r}\frac{B^{T}(\xi)T(u(\xi))}{\pi(z-\xi)}\,d\xi,
\  z\in B_1.
\end{equation}
Since $G(u)(z)=B(G_1(z)+G_2(z))$, we have
\begin{equation}\label{decom3}
|\nabla u|(z)\le C(|G_1(z)|+|G_2(z)|), \ z\in B_1\setminus B_r.
\end{equation}
We can estimate the $L^{2,1}$-norm of $G_1$ in the following Lemma.

\begin{lemma}\label{21-norm-neck2} If $\tau(u)\in L\log L(B_1)$, then
$G_1\in W^{1,1}(B_1)$ and
\begin{equation}\label{21-norm-neck20}
\big\|\nabla G_1\big\|_{L^{1}(B_1)}\le C\big\|\tau(u)\big\|_{L\log L(B_1)},
\end{equation}
and
\begin{equation}\label{21-norm-neck21}
\big\|G_1\big\|_{L^1(B_1)}+\big\|G_1\big\|_{L^{2,\infty}(B_1)}\le C\big\|\tau(u)\big\|_{L^1(B_1)}.
\end{equation}
In particular, $G_1\in L^{2,1}(B_1)$ and the following estimate holds:
\begin{equation}\label{21-norm-neck22}
\big\|G_1\big\|_{L^{2,1}(B_1)}\le C
\big\|\tau(u)\big\|_{L\log L(B_1)}.
\end{equation}
\end{lemma}
\pf It is easy to see from \eqref{g1-form} that 
$$\nabla G_1(z)=\int_{\mathbb C} K(z-\xi) (B^{T}\chi_{B_1\setminus B_r}T(u))(\xi)\,d\xi,$$
where $K:\mathbb C\times \mathbb C\to \mathbb C$ is a kernel of Calderon-Zygmund type. By a result
of Stein \cite{stein}, we know that $\nabla G_1\in L^1(B_1)$ and
\begin{equation}\label{l1-estimate}
\big\|\nabla G_1\big\|_{L^1(B_1)}
\le C\big\|B^{T}\chi_{B_1\setminus B_r}T(u))\big\|_{L\log L(B_1)}
\le C\big\|\tau(u)\big\|_{L\log L(B_1)}.
\end{equation}
This implies \eqref{21-norm-neck20}.
Since $\displaystyle\frac{1}{z}\in L^1\cap L^{2,\infty}(B_1)$, 
 it follows from Young's inequality that 
$G_1\in L^1\cap L^{2,\infty}(B_1)$ and
$$\big\|G_1\big\|_{L^1(B_1)}+\big\|G_1\big\|_{L^{2,\infty}(B_1)}
\le C\big\|B^{T}\chi_{B_1\setminus B_r}T(u)\big\|_{L^1(B_1)}
\le C\big\|\tau(u)\big\|_{L^1(B_1)}.
$$
This gives \eqref{21-norm-neck21}. \eqref{21-norm-neck22} follows from the embedding
$W^{1,1}(B_1)\subset L^{2,1}(B_1)$ (cf. \cite{tartar} and \cite{Helein}):
$$\|f\|_{L^{2,1}(B_1)}\le C\|f\|_{W^{1,1}(B_1)},\ f\in W^{1,1}(B_1),$$
and \eqref{21-norm-neck20} and \eqref{21-norm-neck21}.
 \qed

\medskip
The estimates of $L^{2}$ and $L^{2,1}$-norms of $G_2$ are established in the following Lemma.

\begin{lemma}\label{21-norm-neck1} For any $\lambda>1$ and $0<r\le \frac14$, there exists
a constant $C(\lambda)\le C\lambda^{-1}$, which is monotonically decreasing with respect to $\lambda$, such
that it holds 
\begin{equation}\label{21-norm-neck10}
\big\|G_2\big\|_{L^{2,1}(B_{\lambda^{-1}}\setminus B_{\lambda r})}
\le C|a_{-1}|\big|\ln r\big|+C(\lambda)\big\|G_2\big\|_{L^2(B_1\setminus B_r)},
\end{equation}
\begin{equation}\label{21-norm-neck100}
\big\|\nabla G_2\big\|_{L^{1}(B_{\lambda^{-1}}\setminus B_{\lambda r})}
\le C|a_{-1}|\big|\ln r\big|+C(\lambda)\big\|G_2\big\|_{L^2(B_1\setminus B_r)},
\end{equation}
and
\begin{equation}\label{22-norm-neck10}
\big\|G_2\big\|^2_{L^{2}(B_{\lambda^{-1}}\setminus B_{\lambda r})}
\le C|a_{-1}|^2\big|\ln r\big|+C\lambda^{-2}\big\|G_2\big\|^2_{L^2(B_1\setminus B_r)}.
\end{equation}

\end{lemma}
\pf Direct calculations yield that there exists $C>0$ such that 
$$
\big\|z^{n}\big\|_{L^{2,1}(B_{\lambda^{-1}}\setminus B_{\lambda r})}\le C
\begin{cases} (\lambda r)^{n+1}, 
& {\rm{if}}\ n\le -2, \\
(\frac{1}{\lambda})^{n+1}, & {\rm{if}}\ n\ge 0,\\
|\ln r| & {\rm{if}}\ n=-1.
\end{cases}
$$
This implies that
\begin{eqnarray}\label{21-11}
\big\|G_2\big\|_{L^{2,1}(B_{\lambda^{-1}}\setminus B_{\lambda r})}
&\le& C\Big(|a_{-1}\big|\ln r\big|+\sum_{n\le -2} |a_n|(\lambda r)^{n+1}+\sum_{n\ge 0}|a_n| 
\lambda^{-(n+1)}\Big).
\end{eqnarray}
While
\begin{eqnarray}\label{21-111}
\big\|\nabla G_2\big\|_{L^1(B_{\lambda^{-1}}\setminus B_{\lambda r})}
&\le& C\sum_{n=-\infty}^{\infty} |n||a_n|\int_{B_{\lambda^{-1}}\setminus B_{\lambda r}} |x|^{n-1} \nonumber\\
&\le & C\Big[|a_{-1}|\ln r|+\sum_{n\ge 1} |a_n| \lambda^{-(n+1)}
+\sum_{n\le -2}|a_n|(\lambda r)^{n+1}.
\end{eqnarray}
On the other hand, direct calculations yield
\begin{eqnarray}\label{21-12}
\big\|G_2\big\|_{L^2(B_1\setminus B_r)}^2
&\ge& \pi \Big(\sum_{n\le -2}|a_n|^2\frac{r^{2n+2}-1}{|n+1|}+\sum_{n\ge 0} |a_n|^2\frac{1-r^{2n+2}}{n+1}\Big)\nonumber\\
&\ge& \frac{\pi}2  \Big(\sum_{n\le -2}|a_n|^2\frac{r^{2n+2}}{|n+1|}+\sum_{n\ge 0} |a_n|^2\frac{1}{n+1}\Big).
\end{eqnarray}
It follows from the Cauchy-Schwartz inequality that
\begin{eqnarray}\label{21-13}
\Big(\sum_{n\le -2} |a_n|(\lambda r)^{n+1}+\sum_{n\ge 0}|a_n| 
\lambda^{-(n+1)}\Big)^2
&\le& \Big(\sum_{n\le -2} |a_n|^2\frac{\lambda r^{2(n+1)}}{|n+1|}\Big)
 \Big(\sum_{n\le -2}|n+1|\lambda^{2(n+1)}\Big)\nonumber\\
 &&\ +\Big(\sum_{n\ge 0}|a_n|^2 \frac{1}{n+1}\Big)\Big(\sum_{n\ge 0}
(n+1)\lambda^{-2(n+1)}\Big)\nonumber\\
&\le& C\Big(\sum_{n\le -2}|a_n|^2\frac{r^{2n+2}}{|n+1|}+\sum_{n\ge 0} |a_n|^2\frac{1}{n+1}\Big)\nonumber\\
&\le& C^2(\lambda)\big\|G_1\big\|_{L^2(B_1\setminus B_r)}^2,
\end{eqnarray}
where we have used the fact $\lambda>1$ so that
\begin{equation}\label{c-lambda} \sum_{n\le -2}|n+1|\lambda^{2(n+1)}+\sum_{n\ge 0}
(n+1)\lambda^{-2(n+1)}:=C^2(\lambda)<\infty.
\end{equation}
Note that $C(\lambda)\le C\lambda^{-1}$ and is monotonically decreasing for $\lambda>1$.
It is clear that \eqref{21-norm-neck10} follows by substituting (\ref{21-13}) into (\ref{21-11}) and applying
(\ref{21-12}), and \eqref{21-norm-neck100} follows by substituting (\ref{21-13})
into (\ref{21-111}) and applying (\ref{21-12}). 

It is not hard to check that (\ref{22-norm-neck10}) follows from 
$$\big\|z^{-1}\big\|^2_{L^2(B_{\lambda^{-1}}\setminus B_{\lambda r})}
\le C|\ln r|,$$
and the fact that for any $n\not=-1$,
$$\big\|z^{n}\big\|^2_{L^2(B_{\lambda^{-1}}\setminus B_{\lambda r})}
\le C\lambda^{-2} \big\|z^{n}\big\|^2_{L^2(B_{1}\setminus B_{r})}.
$$
This completes the proof. \qed

\subsection{Hopf differential and its holomorphic approximation}

From (\ref{21-norm-neck10}), in order to control $L^{2}$-norm of $G_2$ in the neck region, we need to
establish a refined estimate of $a_{-1}$, that decays at least as fast as $|\log r|^{-1}$.  This follows
from a delicate holomorphic approximation property of the Hopf differential, which is achieved from an improved version of Wang-Wei-Zhang \cite{WWZ} Proposition 4.1.

Recall that the Hopf differential of a map $u\in H^1(U, N)$ is defined by
$$\mathcal{H}(u):=\big(\frac{\partial u}{\partial z}\big)^2=\big(|\frac{\partial u}{\partial x}|^2
-|\frac{\partial u}{\partial y}|^2\big)-2\langle \frac{\partial u}{\partial x}, \frac{\partial u}{\partial y}\rangle i
\ (=G(u)^T G(u)).$$
It is well-known that if $u$ is harmonic map, then $\mathcal{H}(u)$ is a holomorphic function, i.e.,
$\frac{\partial}{\partial\bar z}(\mathcal{H}(u))=0$.
We let $C^\omega(U)$  denote the space consisting of holomorphic functions on $U$.
Now we will prove the following crucial Lemma on the holomorphic approximation of
$\mathcal H(u)$ for approximate harmonic maps.
\begin{lemma}\label{21-norm-neck3} Assume $u\in W^{1,2}(B_1)$ is an approximate harmonic
map with tension field $\tau(u)\in M^{1,\delta}(B_1)$ for some $1<\delta<2$, satisfying
$\displaystyle\big\|\tau(u)\big\|_{M^{1,\delta}(B_1)}\le 1$ and,
for some $0<r\le\frac14$,
$$\int_{B_1\setminus B_r}|\nabla u|^2\le \epsilon_0^2,$$
where $\epsilon_0>0$ is given by Lemma 2.2. Then\\
(i) the coefficient $a_{-1}$ of the Laurent series 
(2.15) for
$G_2$ in $B_1\setminus B_r$ can be estimated by
\begin{equation}\label{a_{-1}-est}
|a_{-1}|\le C\big(A_0^\frac12(r)+r^\frac12\big),
\end{equation}
where
\begin{equation}\label{A0-form}
A_0(r):=\inf_{h\in C^\omega(B_{2r})}\|\mathcal H(u)-h\|_{L^1(B_{2r})}
+r^{2-\delta}\big\|\tau(u)\big\|_{M^{1,\delta}(B_1)}.
\end{equation}
(ii) there exists $h_1\in C^\omega(B_1)$ such that
\begin{equation}\label{holo-approx}
\big\|\mathcal H(u)-h_1\big\|_{L^1(B_1)}
\le C\Big[A_0(r)|\ln r|+\big(A_0^\frac12(r)+A_0^\frac13(r)\big)
+\big\|\tau(u)\big\|_{M^{1,\delta}(B_1)}\Big].
\end{equation}
(iii) 
for any $\lambda>1$, it holds
\begin{eqnarray}\label{22-norm-neck-bound}
\big\|\nabla u\big\|^2_{L^{2}(B_{\lambda^{-1}}\setminus B_{\lambda r})}
&\le &C(A_0(r)+r)|\ln r|+C\lambda^{-2}\|\nabla u\|^2_{L^2(B_1\setminus B_r)}\nonumber\\
&&+C\big\|\tau(u)\big\|^2_{M^{1,\delta}(B_1)}.
\end{eqnarray}
(iv) if, in addition, $\tau(u)\in L\log L(B_1)$, then for any $\lambda>1$ there holds
\begin{eqnarray}\label{21-norm-neck-bound}
\big\|\nabla u\big\|_{L^{2,1}(B_{\lambda^{-1}}\setminus B_{\lambda r})}
&\le &C(A_0^\frac12(r)+r^\frac12)|\ln r|+C(\lambda)\|\nabla u\|_{L^2(B_1\setminus B_r)}\nonumber\\
&&+C\|\tau(u)\|_{L\log L(B_1)}.
\end{eqnarray}
(v) if, in addition, $\tau(u)\in L\log L(B_1)$, then for any $\lambda>1$ there holds
\begin{eqnarray}\label{W21-norm-neck-bound}
\big\|\nabla^2 u\big\|_{L^{1}(B_{\lambda^{-1}}\setminus B_{\lambda r})}
&\le &C(A_0^\frac12(r)+r^\frac12)|\ln r|+[C(\lambda)+C\lambda^{-1}]\|\nabla u\|_{L^2(B_1\setminus B_r)}\nonumber\\
&&+C\|\tau(u)\|_{L\log L(B_1)}.
\end{eqnarray}
\end{lemma}
\pf Applying Lemma 2.3 with $U=B_1\setminus B_r$ and using the same notations as above, we have
$$B^{T}G(u)(z)=G_1(z)+G_2(z), \ G_2(z)=\sum_{n=-\infty}^\infty a_n z^n, \  z\in B_1\setminus B_r.$$
Since $\mathcal H(u)=G^T(z)G(z)$, we can express
$$\mathcal H(u)=\mathcal H_1(u)+\mathcal H_2(u), \ z\in B_1\setminus B_r,$$
where
$$\mathcal H_1(u)=G^T(z)G_1(z)+G_1^T(z)G_2(z), \ \mathcal H_2(u)=G_2^T(z)G_2(z)
=\sum_{n=-\infty}^n b_n z^n,
\ z\in B_1\setminus B_r,$$
with
$$b_n=\sum_{m=-\infty}^\infty \langle a_m,  a_{n-m}\rangle, \ n\in\mathbb Z.$$
Since $\displaystyle\big\|\tau(u)\big\|_{M^{1,\delta}(B_1)}\le 1$,  it follows from  \eqref{G1-morrey-est1} that
\begin{eqnarray}\label{holo-approx1}
\big\|\mathcal H(u)-\mathcal H_2(u)\big\|_{L^1(B_1\setminus B_r)}
&\le & \big\|G_1\big\|_{L^{2}(B_1\setminus B_r)}
\Big(2\big\|G\big\|_{L^{2}(B_1\setminus B_r)}+\big\|G_1\big\|_{L^{2}(B_1\setminus B_r)}\Big)\nonumber\\
&\le & C\big\|\tau(u)\big\|_{M^{1,\delta}(B_1)}\Big(\epsilon_0+
\big\|\tau(u)\big\|_{M^{1,\delta}(B_1)}\Big)\nonumber\\
&\le & C\big\|\tau(u)\big\|_{M^{1,\delta}(B_1)}.
\end{eqnarray}
For simplicity, denote $h_2=\mathcal H_2(u)\in C^\omega(B_1\setminus B_r)$.
Set $$\displaystyle h_1(z)=\sum_{n=0}^\infty b_n z^n\in C^\omega(B_1).$$
For an arbitrary small $\epsilon>0$, choose $\displaystyle h_{0}=\sum_{n=0}^\infty c_n z^n\in C^\omega(B_{2r})$ such that
$$\|\mathcal H(u)-h_0\|_{L^1(B_{2r})}\le \inf_{h\in C^\omega(B_{2r})}\|\mathcal H(u)-h\|_{L^1(B_{2r})}+\epsilon.$$
As in \cite{WWZ} Lemma 4.3, set $r_1=\frac54 r, r_2=\frac32 r, r_3=\frac74 r$ and estimate
\begin{eqnarray}\label{holo-approx2}
\big\|\mathcal H(u)-h_1\big\|_{L^1(B_1)}
&\le& \big\|\mathcal H(u)-h_1\big\|_{L^1(B_1\setminus B_{r_2})}
+\big\|\mathcal H(u)-h_1\big\|_{L^1(B_{r_2})}\nonumber\\
&\le& \big\|\mathcal H(u)-h_2\big\|_{L^1(B_1\setminus B_{r_2})}
+\big\|\mathcal H(u)-h_0\big\|_{L^1(B_{r_2})}\nonumber\\
&+&\big\|h_2-h_1\big\|_{L^1(B_1\setminus B_{r_2})}
+\big\|h_0-h_1\big\|_{L^1(B_{r_2})}\nonumber\\
&\le& C\big\|\tau(u)\big\|_{M^{1,\delta}(B_1)}+\inf_{h\in C^\omega(B_{2r})}\|\mathcal H(u)-h\|_{L^1(B_{2r})}+\epsilon\nonumber\\
&+&C\big[|b_{-1}|+|b_{-2}||\ln r_2|+\sum_{n\ge 3}|b_{-n}|r_2^{2-n}+\sum_{n\ge 0} |b_n-c_n|r_2^{n+2}\big].
\end{eqnarray}
Observe that, by an argument similar to (\ref{holo-approx1}) and Lemma 2.3, we have
\begin{eqnarray} \|\mathcal H(u)-h_2\|_{L^1(B_{2r})}&\le& C\|G_1\|_{L^2(B_{2r})}
\big(\|G\|_{L^{2}(B_{2r})}+\|G_1\|_{L^2(B_{2r})}\big)\nonumber\\
&\le& Cr^{2-\delta}\big\|G_1\big\|_{M^{2,2\delta-2}(B_1)}
\Big[\|G\|_{L^{2}(B_{2r})}+r^{2-\delta}\big\|G_1\big\|_{M^{2,2\delta-2}(B_1)}\Big]\nonumber\\
&\le& Cr^{2-\delta}\big\|\tau(u)\big\|_{M^{1,\delta}(B_1)}
\Big(\|\nabla u\|_{L^{2}(B_1)}+r^{2-\delta}\big\|\tau(u)\big\|_{M^{1,\delta}(B_1)}\Big)\nonumber\\
&\le & Cr^{2-\delta}\big\|\tau(u)\big\|_{M^{1,\delta}(B_1)}.\label{holo-approx12}
\end{eqnarray}
Hence we obtain
$$\|h_2-h_0\|_{L^1(B_{2r}\setminus B_r)}
\le \|\mathcal H(u)-h_2\|_{L^1(B_{2r})}+\|\mathcal H(u)-h_0\|_{L^1(B_{2r})}
\le CA_0(r)+\epsilon.
$$
Applying Fubini's theorem, we may assume that
$$r_j\int_{|z|=r_j}|h_2-h_0||dz|\le C\|h_2-h_0\|_{L^1(B_{2r}\setminus B_r)}\le CA_0(r)+\epsilon,
\ j=1, 2, 3.$$
By the Laurent series coefficient formula, we can then estimate
\begin{eqnarray*}
\big|b_{-n}\big|&=&\frac{1}{2\pi}
\big|\int_{|z|=r_1}z^{n-1}h_2(z)\,dz\big|
=\frac{1}{2\pi}\big|\int_{|z|=r_1}z^{n-1}(h_2(z)-h_0(z))\,dz\big|\\
&\le& Cr_1^{n-2} r_1\int_{|z|=r_1}|h_2-h_0||dz|\le (CA_0(r)+\epsilon)r_1^{n-2},
\end{eqnarray*}
for all $n\ge 1$, and for $n\ge 0$
\begin{eqnarray*}
|b_n-c_n|&=&\frac{1}{2\pi}\big|\int_{|z|=r_3}z^{-(n+1)}h_2(z)\,dz\big|
=\frac{1}{2\pi}\big|\int_{|z|=r_3}z^{-(n+1)}(h_2-h_0)(z)\,dz\big|\\
&\le& Cr_3^{-(n+2)} r_3\int_{|z|=r_3}|h_2-h_0||dz|\le (CA_0(r)+\epsilon)r_3^{-(n+2)}.
\end{eqnarray*}
Putting all these estimate into \eqref{holo-approx2} and sending $\epsilon$ to $0$, we obtain
\begin{equation}\label{holo-approx3}
\big\|\mathcal H(u)-h_1\big\|_{L^1(B_1)}
\le C\Big(\big\|\tau(u)\big\|_{M^{1,\delta}(B_1)}+A_0(r)|\ln r|+|b_{-1}|\Big).
\end{equation}
The estimate of $a_{-1}$ can be done as in \cite{WWZ} Lemma 4.3 with mild modifications. 
Denote 
$$Q(z)=(e_1(u(z)), \cdots, e_l(u(z))),\ \   q(z)=\frac{1}{2\pi |z|}\int_{|w|=|z|}Q(w)\,dw.$$
Then 
$$G_2=B^T G-G_1=B^T Q^T\frac{\partial u}{\partial z}-G_1
=(B-\mathbb I_l)^TQ^T\frac{\partial u}{\partial z}+q^T\frac{\partial u}{\partial z}+(Q-q)^T\frac{\partial u}{\partial z}-G_1.$$
Thus we have that, for
$\rho=\sqrt{\frac{r}2}$, 
\begin{eqnarray*}
(2\pi \ln 2) a_{-1}&=&\int_{B_{2\rho}\setminus B_\rho}\frac{G_2(z)}{\bar z}\, \frac{i}2dz\wedge d\bar z\\
&=&\int_{B_{2\rho}\setminus B_\rho}\frac{q^T}{\bar z}\frac{\partial u}{\partial z}\, \frac{i}2dz\wedge d\bar z\\
&+&\int_{B_{2\rho}\setminus B_\rho}\big[(B-\mathbb I_l)^T\frac{Q^T}{\bar z}\frac{\partial u}{\partial z}
+\frac{(Q-q)^T}{\bar z}\frac{\partial u}{\partial z}-\frac{G_1}{\bar z}\big]\, \frac{i}2dz\wedge d\bar z\\
&=& I_1+I_2.
\end{eqnarray*}
As pointed out by \cite{WWZ} Lemma 4.3,  each component of $I_1\in \mathbb C^l$  is a real number. Hence
\begin{eqnarray*}
&&2\pi \ln 2 |{\rm{Im}}a_{-1}|=|I_2|\\
&\le& \int_{B_{2\rho}\setminus B_\rho}\big[|(B-\mathbb I_l)^T\frac{Q^T}{\bar z}\frac{\partial u}{\partial z}|
+|\frac{(Q-q)^T}{\bar z}\frac{\partial u}{\partial z}|+|\frac{G_1}{\bar z}|\big]|dz\wedge d\bar z|\\
&\le& C\Big[\|B-\mathbb I_l\|_{L^\infty(B_{2\rho}\setminus B_\rho)}
\|\frac{1}{\bar z}\|_{L^2(B_{2\rho}\setminus B_\rho)}\|\frac{\partial u}{\partial z}\|_{L^2(B_{2\rho}\setminus B_\rho)}\\
&&+\|\frac{1}{\bar z}\|_{L^2(B_{2\rho}\setminus B_\rho)}\|G_1\|_{L^2(B_{2\rho}\setminus B_\rho)}
+\big\|\frac{(Q-q)^T}{\bar z}\big\|_{L^2(B_{2\rho}\setminus B_\rho)}\|\frac{\partial u}{\partial z}\|_{L^2(B_{2\rho}\setminus B_\rho)}\Big]\\
&\le& C\epsilon_0\|G\|_{L^2(B_{2\rho}\setminus B_\rho)}+C\|G_1\|_{L^2(B_{2\rho}\setminus B_\rho)}\\
&\le& C\epsilon_0\|G_2\|_{L^2(B_{2\rho}\setminus B_\rho)}+C\|G_1\|_{L^2(B_{2\rho}\setminus B_\rho)},
\end{eqnarray*}
where we have used the inequality from Lemma 2.2:
$$\|B-\mathbb I_l\|_{L^\infty(B_{2\rho}\setminus B_\rho)}\le C\|\nabla u\|_{L^2(B_1\setminus B_r)}\le C\epsilon_0,$$
and the Hardy inequality:
$$
\Big\|\frac{(Q-q)^T}{\bar z}\Big\|_{L^2(B_{2\rho}\setminus B_\rho)}
\le C\|\nabla Q\|_{L^2(B_{2\rho}\setminus B_\rho)}\le C\|\nabla u\|_{L^2(B_1\setminus B_r)}\le
C\epsilon_0.$$
As shown by \cite{WWZ} Lemma 4.3, it holds
$$\|G_2\|_{L^2(B_{2\rho}\setminus B_\rho)}^2
\le 2\pi \ln 2|a_{-1}|^2+C\rho^2\|G_2\|_{L^2(B_1\setminus B_r)}^2
\le 2\pi \ln 2|a_{-1}|^2+Cr,$$
where we have used the fact that
\begin{eqnarray*}\|G_2\|_{L^2(B_1\setminus B_r)}^2&\le& 2\big(\|G\|_{L^2(B_1\setminus B_r)}^2
+\|G_1\|_{L^2(B_1\setminus B_r)}^2\big)\\
&\le& C\Big(\|\nabla u\|_{L^2(B_1\setminus B_r)}^2
+\big\|\tau(u)\big\|^2_{M^{1,\delta}(B_1)}\Big)\le C.
\end{eqnarray*}
Therefore we obtain
\begin{equation}\label{a_{-1}-estimate1}
 |{\rm{Im}}a_{-1}|\le  C\epsilon_0 \big(2\pi \ln 2|a_{-1}|^2+Cr\big)^\frac12+C
\|G_1\|_{L^2(B_{2\rho}\setminus B_\rho)}.
\end{equation}
On the other hand, it was proven by \cite{WWZ} Lemma 4.3 that
$$\sum_{n\ge 0} |a_n||a_{-2-n}|\le Cr\|G_2\|_{L^2(B_1\setminus B_r)}^2.$$
Since
$$|\langle a_{-1}, a_{-1}\rangle|
\le |b_{-2}|+2\sum_{n\ge 0} |a_n|a_{-2-n}|,$$
it then follows
$$|\langle a_{-1}, a_{-1}\rangle|\le CA_0(r)+Cr\|G_2\|_{L^2(B_1\setminus B_r)}^2
\le C(A_0(r)+r).$$
Since $|a_{-1}|^2={\rm{Re}}\langle a_{-1}, a_{-1}\rangle+2|{\rm{Im}}a_{-1}|^2$,
we obtain
$$
|a_{-1}|^2\le C(A_0(r)+r)+C\epsilon_0^2\big(2\pi \ln 2|a_{-1}|^2+Cr\big)+C
\|G_1\|_{L^2(B_{2\rho}\setminus B_\rho)}^2,
$$
and hence
$$|a_{-1}|^2\le C(A_0(r)+r)+C
\|G_1\|_{L^2(B_{2\rho}\setminus B_\rho)}^2,$$
provided $\epsilon_0>0$ is chosen to be sufficiently small. 

It follows from Lemma 2.3  and \eqref{G1-morrey-est1} that 
\begin{eqnarray*}
\|G_1\|_{L^2(B_{2\rho}\setminus B_\rho)}^2
&\le& \rho^{4-2\delta}\big\|G_1\big\|^2_{M^{2,2\delta-2}(B_1)}
\le C\rho^{4-2\delta} \big\|\tau(u)\big\|^2_{M^{1,\delta}(B_1)}\\
&\le& Cr^{2-\delta}\big\|\tau(u)\big\|_{M^{1,\delta}(B_1)}^2\le CA_0(r).
\end{eqnarray*}
Thus we finally obtain
$$|a_{-1}|^2\le C(A_0(r)+r),$$
which yields (\ref{a_{-1}-est}).

The estimate of $b_{-1}$ above needs to be refined as follows. First, note that
$$
|b_{-1}|\le 2\big(|a_0||a_{-1}|+\sum_{n\ge 0}|a_n||a_{-1-n}|\big).$$
Recall that it was proved by \cite{WWZ} Lemma 4.3 that
$$
\sum_{n\ge 0}|a_n||a_{-1-n}|\le Cr\|G_2\|_{L^2(B_1\setminus B_r)}^2\le Cr.$$
Hence
$$|b_{-1}|\le C(|a_{-1}|+r)\le C\big(A_0(r)^\frac12+r^\frac12\big).$$
Thus
$$|b_{-1}|\le C\min\big\{\frac{A_0(r)}{r}, A_0(r)^\frac12+r^\frac12\big\}
\le C\big(A_0(r)^\frac12+A_0(r)^\frac13\big).$$
Putting this estimate into (\ref{holo-approx2}) yields (\ref{holo-approx}).

It is not hard to check that \eqref{22-norm-neck-bound} follows from  the estimate (\ref{22-norm-neck10}) for $G_2$, the estimate \eqref{G1-morrey-est1} for $G_1$, \eqref{decom3}
and the estimate (\ref{a_{-1}-est}) for $a_{-1}$. While \eqref{21-norm-neck-bound} follows from
(\ref{21-norm-neck10}) for $G_2$, \eqref{21-norm-neck22} for $G_1$,  \eqref{decom3}
and the estimate (\ref{a_{-1}-est}) for $a_{-1}$.

To prove (\ref{W21-norm-neck-bound}), we first observe that differentiating \eqref{decom2} yields
\begin{equation}\label{decom4}
|\nabla^2 u|\le C|\nabla B|(|G_1|+|G_2|)+C(|\nabla G_1|+|\nabla G_2|)+C\sum_{\alpha=1}^l
(|G_1|+|G_2|)|\nabla e_\alpha|,  \ {\rm{in}}\ B_1\setminus B_r.
\end{equation}
Note that by \eqref{col-frame} and Lemma \ref{conjugate}, we have
$$\big\|\nabla B\big\|_{L^2(B_1\setminus B_r)}
+\sum_{\alpha=1}^l\big\|\nabla e_\alpha\big\|_{L^2(B_1\setminus B_r)}
\le C\big\|\nabla u\big\|_{L^2(B_1\setminus B_r)}\le C\epsilon_0.$$
Thus by H\"older's inequality we obtain
\begin{eqnarray*}
\int_{B_{\lambda^{-1}}\setminus B_{\lambda r}}|\nabla^2 u|
&\le& C\Big(\big\|\nabla G_1\big\|_{L^1(B_{\lambda^{-1}}\setminus B_{\lambda r})}
+\big\|\nabla G_2\big\|_{L^1(B_{\lambda^{-1}}\setminus B_{\lambda r})}\Big)\\
&&+C\Big(\big\|G_1\big\|_{L^2(B_{\lambda^{-1}}\setminus B_{\lambda r})}
+\big\|G_2\big\|_{L^2(B_{\lambda^{-1}}\setminus B_{\lambda r})}\Big).
\end{eqnarray*}
Hence applying the estimates \eqref{21-norm-neck20} for $\nabla G_1$ and \eqref{21-norm-neck21} for $G_1$,
and the estimates \eqref{21-norm-neck100} $\nabla G_2$ and \eqref{22-norm-neck10} for $G_2$,
we obtain (\ref{W21-norm-neck-bound}). The proof is now complete. \qed

\begin{theorem}\label{holo-approx5} For an approximate harmonic map
$u\in W^{1,2}(B_1, N)$ with tension field $\tau(u)\in M^{1,\delta}(B_1)$
for some $1<\delta<2$, assume $E(u,B_1)\le m\epsilon_0^2$ for some positive integer $m\ge 1$,
here $\epsilon_0$ is given by Lemma 2.2.
Then there exists $C_m>0$ such that for any $0<r\le 1$, there is a holomorphic function 
$h\in C^\omega(B_\frac{r}4)$ such that
\begin{equation}\label{holo-approx6}
\big\|\mathcal H(u)-h\big\|_{L^1(B_\frac{r}4)}\le C_m 
\Big(r^{2-\delta}\big\|\tau(u)\big\|_{M^{1,\delta}(B_r)}\Big)^{3^{1-m}}.
\end{equation}
\end{theorem}
\pf It follows from \cite{WWZ} Proposition 4.1 with suitable modifications. Here we only sketch it.
By simple scaling argument, it suffices to prove \eqref{holo-approx6} for $r=1$.

Note that \eqref{holo-approx6} trivially holds, if $\displaystyle\big\|\tau(u)\big\|_{M^{1,\delta}(B_1)}\ge 1$.
Thus we may assume
\begin{equation}\label{small-phi}
\big\|\tau(u)\big\|_{M^{1,\delta}(B_1)}\le 1,
\end{equation}
and prove \eqref{holo-approx6} by an induction on $m$. If $m=1$, then as in the proof of (2.30), with
the domain $B_1\setminus B_r$ replaced by $B_1$, we conclude that there exists $h\in C^\omega(B_1)$ such that
$$\big\|\mathcal H(u)-h\big\|_{L^1(B_1)}\le C \big\|\tau(u)\big\|_{M^{1,\delta}(B_1)}.$$
This yields (\ref{holo-approx6}).
For $m\ge 2$, assume \eqref{holo-approx6} holds for any approximate harmonic map $u$ satisfying
(\ref{small-phi}) and $E(u, B_1)\le (k-1)\epsilon_0^2$ for $k\ge 2$, we want to prove (\ref{holo-approx6}) for $m=k$.
For this, we may further assume 
\begin{equation}\label{small-energy}
E(u, B_\frac14)> \epsilon_0^2.
\end{equation}
For, otherwise, $v(x)=u(\frac{x}4): B_1\to N$ is an approximate harmonic map, with tension field
$\tau(v)(x)=\frac{1}{16}\tau(u)(\frac{x}4)$, satisfying
$$E(v, B_1)=E(u, B_\frac14)\le \epsilon_0^2;\  \ \big\|\tau(v)\big\|_{M^{1,\delta}(B_1)}
\le 4^{{\delta-2}}\big\|\tau(u)\big\|_{M^{1,\delta}(B_\frac14)}\le 1.
$$
Hence, from the case $m=1$, there exists $\hat{h}\in C^\omega(B_1)$ such that
$$\big\|\mathcal H(v)-\hat{h}\big\|_{L^1(B_1)}
\le C\big\|\tau(v)\big\|_{M^{1,\delta}(B_1)}\le C\big\|\tau(u)\big\|_{M^{1,\delta}(B_1)}.$$
By scalings, this implies that (\ref{holo-approx6}) holds for $\mathcal{H}(u)$ and
$h(x)=16\hat{h}(4x)\in C^\omega(B_\frac14)$,
since
$$\big\|\mathcal H(u)-h\big\|_{L^1(B_\frac14)}=\big\|\mathcal H(v)-\hat{h}\big\|_{L^1(B_1)}.$$
Define the energy concentration function of $u$ by
$$\mathcal E(r)=\sup_{B_r(x)\subset B_1} E(u, B_r(x)),$$
and divide the proof into two cases.\\
(i) $\mathcal E(\frac18)\le (k-1)\epsilon_0^2$: For any $x_0\in B_\frac12$, consider $w(x)=u(x_0+\frac18 x):
B_1\to N$. Then it is easy to check that $w$ is an approximate harmonic map, with tension field
$\tau(w)(x)=\frac{1}{64}\tau(u)(\frac{x}8)$, satisfying
$$\begin{cases}E(w, B_1)=E(u, B_\frac18(x_0))\le \mathcal E(\frac18)\le (k-1)\epsilon_0^2; \\
\big\|\tau(w)\big\|_{M^{1,\delta}(B_1)}\le 8^{\delta-2}\big\|\tau(u)\big\|_{M^{1,\delta}(B_\frac18(x_0))}
\le \big\|\tau(u)\big\|_{M^{1,\delta}(B_1)}\le 1.
\end{cases}
$$
Thus, by the induction hypothesis there exists $\tilde{h}\in C^\omega(B_\frac14)$ such that
$$
\big\|\mathcal H(w)-\tilde{h}\big\|_{L^1(B_\frac14)}\le C_{k-1}\big\|\tau(w)\big\|_{M^{1,\delta}(B_1)}^{3^{2-k}}
\le C_{k-1}\big\|\tau(u)\big\|_{M^{1,\delta}(B_1)}^{3^{2-k}}.
$$
By scalings, this implies that for any $x_0\in B_\frac12$, there exists $\tilde{h}_{x_0}(z)
=64\tilde{h}(8(z-x_0))\in C^\omega(B_\frac{1}{32}(x_0))$ such that
$$
\big\|\mathcal H(u)-\tilde{h}_{x_0}\big\|_{L^1(B_\frac1{32}(x_0))}
\le C_{k-1}\big\|\tau(u)\big\|_{M^{1,\delta}(B_1)}^{3^{2-k}}.
$$
Now, applying the gluing Lamma 4.4 of \cite{WWZ}, we conclude that there exists $h\in C^\omega(B_\frac14)$ such that 
$$
\big\|\mathcal H(u)-h\big\|_{L^1(B_\frac14)}
\le C\big(\frac{1/4}{1/32}\big)^3 C_{k-1}\big\|\tau(u)\big\|_{M^{1,\delta}(B_1)}^{3^{2-k}}
\le 8^3 CC_{k-1}\big\|\tau(u)\big\|_{M^{1,\delta}(B_1)}^{3^{2-k}}.
$$
This proves (\ref{holo-approx6}) for $m=k$. Now we proceed with the difficult case.\\
\noindent(ii) $\mathcal E(\frac18)> (k-1)\epsilon_0^2$: There exists $0<r_1<\frac18$ and $B_{r_1}(x_1)\subset B_1$
such that
$$E(u, B_{r_1}(x_1))=\mathcal E(r_1)=(k-1)\epsilon_0^2.$$
Since $E(u, B_1)\le k\epsilon_0^2$ and $E(u, B_\frac14)>\epsilon_0^2$, this implies
that $B_{r_1}(x_1)\cap B_\frac14\not=\emptyset$. Hence $|x_1|\le \frac14+r_1$, and $B_{4r_1}(x_1)\subset
B_1$. For any $z_0\in B_{3r_1}(x_1)$, define
$\widehat{w}(z)=u(z_0+r_1 z): B_1\to N$. Then $\widehat{w}$ is an approximate harmonic map, with
tension field $\tau(\widehat{w})(z)=r_1^2\tau(u)(z_0+r_1z)$, satisfying
$$\begin{cases}
E(\widehat{w}, B_1)=E(u, B_{r_1}(z_0))\le \mathcal E(r_1)=(k-1)\epsilon_0^2,\\
\big\|\tau(\widehat{w})\big\|_{M^{1,\delta}(B_1)}
\le r_1^{2-\delta} \big\|\tau({u})\big\|_{M^{1,\delta}(B_1)}\ (\le 1).
\end{cases}
$$
Now we can repeat the same argument as in the case (i), with $w$ replaced by $\widehat{w}$, to conclude
that there exists a holomorphic function ${h}_{x_1}\in C^\omega(B_{2r_1}(x_1))$
such that
\begin{equation}\label{holo-approx7}
\big\|\mathcal H(u)-{h}_{x_1}\big\|_{L^1(B_{2r_1}(x_1))}
\le 8^3 CC_{k-1} \Big(r_1^{2-\delta} \big\|\tau({u})\big\|_{M^{1,\delta}(B_1)}\Big)^{3^{2-k}}.
\end{equation}
Now we rescale $u$ again. Set $\rho=\frac12+r_1$ and $r=\frac{r_1}{\rho}$,
and define $\widetilde{w}(z)=u(x_1+\rho z)$. Then $\tau(\widetilde{w})(z)=\rho^2\tau(u)(x_1+\rho z)$, and
$$
\begin{cases}
E(\widetilde{w}, B_1\setminus B_{r})=E(u, B_{\rho}(x_1)\setminus B_{r_1}(x_1))
\le E(u, B_1)-E(u, B_{r_1}(x_1))\le\epsilon_0^2,\\
\big\|\tau(\widetilde{w})\big\|_{M^{1,\delta}(B_1)}\le
\rho^{2-\delta}\big\|\tau(u)\big\|_{M^{1,\delta}(B_1)}\ (\le 1).
\end{cases}
$$
Let $\widetilde{A}_0(r)$ be the quantity, associated with $\widetilde{w}$ and $\tau(\widetilde{w})$,
given by (2.27). 
Define $\widetilde{h}(z)=\rho^2h_{x_1}(x_1+\rho z)\in C^\omega(B_{2r})$. Then by scalings and
(\ref{holo-approx7}) we have
\begin{eqnarray}\label{holo-approx8}
\widetilde{A}_0(r)&\le&\big\|\mathcal H(\widehat{w})-\widetilde{h}\big\|_{L^1(B_{2r})}+r^{2-\delta}
\big\|\tau(\widetilde{w})\big\|_{M^{1,\delta}(B_1)}\nonumber\\
&\le& \big\|\mathcal H(u)-{h}_{x_1}\big\|_{L^1(B_{2r_1}(x_1))}+r_1^{2-\delta}
\big\|\tau(u)\big\|_{M^{1,\delta}(B_1)}\nonumber\\
&\le& CC_{k-1} \Big(r_1^{2-\delta}
\big\|\tau(u)\big\|_{M^{1,\delta}(B_1)}\Big)^{3^{2-k}}.
\end{eqnarray}
Applying Lemma 2.6 (\ref{holo-approx}) to $\widetilde{w}$, we conclude that there exists
$\widetilde{h}_1\in C^\omega(B_1)$ such that
\begin{eqnarray}\label{holo-approx9}
\big\|\mathcal H(\widetilde{w})-\widetilde{h}_1\big\|_{L^1(B_1)}
&\le& C\Big[\widetilde{A}_0(r)|\ln r|+\widetilde{A}_0^\frac12(r)+\widetilde{A}_0^\frac13(r)
+\big\|\tau(u)\big\|_{M^{1,\delta}(B_1)}\Big]\nonumber\\
&\le& C\Big[ CC_{k-1}r_1^{\frac{2-\delta}{3^k}} |\ln r| \big\|\tau(u)\big\|_{M^{1,\delta}(B_1)}^{3^{2-k}}
+C\big\|\tau(u)\big\|_{M^{1,\delta}(B_1)}^{3^{1-k}}\Big]\nonumber\\
&\le& C_k \big\|\tau(u)\big\|_{M^{1,\delta}(B_1)}^{3^{1-k}},
\end{eqnarray}
where we have used the fact that $\frac12\le\rho\le 1$ and $0<r<\frac14$, and hence
$$r_1^{\frac{2-\delta}{3^k}} |\ln r|
= \rho^{\frac{2-\delta}{3^k}}r^{\frac{2-\delta}{3^k}}|\ln r|
\le r^{\frac{2-\delta}{3^k}}|\ln r|\le C.
$$
Define $h_1(z)=\rho^{-2}\widetilde{h}_1(\rho^{-1}(z-x_1))\in C^\omega(B_{\rho}(x_1))$. 
It immediately follows from $B_\frac14\subset B_{\rho}(x_1)$
and \eqref{holo-approx9} that
$$\big\|\mathcal H(u)-{h}_1\big\|_{L^1(B_{\frac14})}\le
\big\|\mathcal H(u)-{h}_1\big\|_{L^1(B_{\rho}(x_1))}
=\big\|\mathcal H(\widetilde{w})-\widetilde{h}_1\big\|_{L^1(B_1)}
\le C_k \big\|\tau(u)\big\|_{M^{1,\delta}(B_1)}^{3^{1-k}}.
$$
This completes the proof. \qed
\section{Estimates over the neck region and proof of Theorem \ref{bubbling1}}
\setcounter{equation}{1}
\setcounter{theorem}{2}
This section is devoted to the proof of the energy identity and the bubble tree convergence
results as stated in Theorem \ref{bubbling1}.

\noindent{\bf Proof of Theorem 1.2}:

Since $M^{1,\delta_1}(\Sigma)\subset M^{1,\delta_2}(\Sigma)$ for any $0\le \delta_1\le\delta_2<2$,
we may, from now on, simply assume the condition (\ref{tension-bound}) holds for some $1<\delta<2$.
Let $\epsilon_0>0$ be given by Lemma \ref{epsilon-regularity1} and define the concentration set
$$\mathcal S=\bigcap_{r>0}\Big\{x\in\Sigma: \ \liminf_{n}\int_{B_r(x)}|\nabla u_n|^2>\epsilon_0^2\Big\}.$$
Then it is well-known (cf. e.g., \cite{Lin-Wang2}) that there exists $m_0\le \frac{E_0}{\epsilon_0^2}$ such that
$\Sigma$ is a finite set of at most $m_0$ points, here $\displaystyle E_0=\sup_{n} \int_\Sigma |\nabla u_n|^2\,dv_g$.
Moreover, it follows from Lemma \ref{epsilon-regularity1} that for any compact set $K\subset\Sigma\setminus\mathcal S$,
\begin{equation}\label{uniform-est}
\sup_n\Big(\big[u_n\big]_{C^{2-\delta}(K)}+\big\|u_n\big\|_{W^{1,q}(K)}\Big)\le 
C\big(q, \delta, \epsilon_0, E_0, {\rm{dist}}(K, \partial({\Sigma\setminus\mathcal S}))\big),
\end{equation}
holds for all $2\le q<\frac{\delta}{\delta-1}$. Thus we may assume that
\begin{equation}\label{strong-convergence}
u_n\rightarrow u \ \ {\rm{in}}\  \ C^{2-\delta}_{\rm{loc}}\cap W^{1,q}_{\rm{loc}}(\Sigma\setminus\mathcal S),
\ \forall\ q\in [1, \frac{\delta}{\delta-1}).
\end{equation}
It is clear that $u$ is an approximate harmonic map with tension field $\tau\in M^{1,\delta}(\Sigma)$,
Applying Lemma \ref{epsilon-regularity1} to $u$, one concludes that $u\in 
C^{2-\delta}\cap W^{1,q}(\Sigma)$ for any $1\le q<\frac{\delta}{\delta-1}$. This proves (i). 

Set $\mathcal S=\big\{x_1,\cdots, x_{m_0}\big\}$ and 
$r(\mathcal S)=\frac12\min\big\{|x-y|: \ x, y\in \mathcal S, \ x\not=y\big\}>0$.
Then for $x_1\in\mathcal S$, for some large constant $C_0>1$ (to be determined later)
there exist $0<r_n^1\le r(\mathcal S)$ and $x_n^1\in B_{r(\mathcal S)}(x_1)$ such that
\begin{equation}\label{exhaust1}
\int_{B_{r_n^1}(x_n^1)}|\nabla u_n|^2=\mathcal E(u_n, r_n^1)=\max\Big\{\int_{B_{r_n^1}(x)}|\nabla u_n|^2: 
\ B_{r_n^1}(x)\subset B_{r(\mathcal S)}(x_1)\Big\}
=\frac{\epsilon_0^2}{C_0}.
\end{equation} 
It is readily seen that $x_n^1\rightarrow x_1$ and $r_n^1\rightarrow 0$. Define the blowing up sequence
$$v_n^1(x)=u_n(x_n^1+r_n^1 x): \Omega_n\equiv (r_n^1)^{-1}\big(B_{r(\mathcal S)}(x_1)\setminus\{x_n^1\}\big)\to
N.$$ 
Then $v_n^1$ is an approximate harmonic map, with tension field
 $\tau(v_n^1)(x)=(r_n^1)^2\tau(u_n)(x_n^1+r_n^1 x)$, such that\\
 (i) $\displaystyle\big\|\tau(v_n^1)\big\|_{M^{1,\delta}(\Omega_n)}\le (r_n^1)^{2-\delta}\big\|\tau(u_n)\big\|_{M^{1,\delta}(\Sigma)}
 \rightarrow 0$, and $\Omega_n\rightarrow\R^2$.\\
 (ii) $\displaystyle\int_{\Omega_n}|\nabla v_n^1|^2\le \int_{\Sigma}|\nabla u_n|^2\le E_0.$\\
 (iii) 
 $\displaystyle\int_{B_1(x)}|\nabla v_n^1|^2\le \int_{B_1(0)}|\nabla v_n^1|^2=\frac{\epsilon_0^2}{C_0}, \ \forall\ x\in\Omega_n.$ In particular, if we choose $C_0$ sufficiently large, then it holds
 $$\int_{B_{10}(x)\cap\Omega_n}|\nabla v_n^1|^2\le \epsilon_0^2, \ \forall \ x\in\Omega_n.$$
Applying Lemma \ref{epsilon-regularity1}, we conclude that there exists a nontrivial harmonic map
$\omega_1\in W^{1,2}(\R^2,N)$ such that $v_n^1\rightarrow \omega_1$ in
$C^{2-\delta}_{\rm{loc}}(\R^2)\cap W^{1,q}_{\rm{loc}}(\R^2)$ for any $1<q<\frac{\delta}{\delta-1}$.
It is well-known \cite{SaU} that $\omega_1$ can be lifted into a nontrivial harmonic map from $\mathbb S^2$ to
$N$. Repeating this process up to finitely many times, similar to \cite{BC} and \cite{Qing},
we can find all possible bubbles near $\mathcal S$, $\{\omega_i\}$ for $1\le i\le m$, and 
sequences of blowing up points and scales $\{x_n^i\}\subset\Sigma$ and $\{r_n^i\}\subset\R_+$
for $1\le i\le m$ such that 
\begin{equation}\label{bubble-convergence}
v_n^i(\cdot)=u_n(x_n^i+r_n^i \cdot)\rightarrow \omega_i \ {\rm{in}}\ C^{2-\delta}_{\rm{loc}}(\R^2)\cap W^{1,q}_{\rm{loc}}(\R^2), \ \forall\ 1<q<\frac{\delta}{\delta-1},
\end{equation}
for all $1\le i\le m$. Moreover, the blowing up points and scales satisfy the relation \eqref{scale-separation}.

To prove the energy identity \eqref{bubbling2} or \eqref{bubble2.0},  it suffices to show that 
there is no energy concentration in the neck regions, which is either the shrinking annulus 
$A_n=B_{r_n^2}(x_n^2)\setminus B_{r_n^1}(x_n^1)$ between
two bubbles with almost same bubbling point but different bubbling scales, or the 
annulus between the body region and bubbling region $A_n=B_{\eta}(x_n^1)\setminus B_{r_n^1}(x_n^1)$.
We proceed the proof by dividing it into two cases.\\
(a) The neck region $A_n$ between two bubbles formed at the same center, but with
different blowing up scales, i.e., $A_n=B_{r_n^2}(x_n^2)\setminus B_{r_n^1}(x_n^1)$, with $\displaystyle
\mu_n=\frac{r_n^2}{r_n^1}\rightarrow\infty$ and $\displaystyle\frac{|x_n^2-x_n^1|}{r_n^1+r_n^2}=0$. For simplicity, assume $x_n^1=x_n^2=x_1$. In this case, we have
$$E(u_n, A_n)=\int_{A_n}|\nabla u_n|^2\le \epsilon_0^2.$$
Denote $\displaystyle E_0=\sup_{n}E(u_n, \Sigma)$ and $\displaystyle m=\big[\frac{E_0}{\epsilon_0}\big]$.
Set $r_n=\mu_n^{-1}=\frac{r_n^1}{r_n^2}\rightarrow 0$, and 
define the blowing up sequence $v_n(x)=u_n(x_1+r_n^2 x): B_1\to N$. Then $v_n$ is a sequence of approximate harmonic maps,
with tension fields $\tau(v_n)(x)=(r_n^2)^2\tau(u_n)(x_1+r_n^2 x)$, satisfying
$$E(v_n, B_1)\le E(u_n, \Sigma)\le m\epsilon_0^2;\ \  E(v_n, B_1\setminus B_{r_n})
=E(u_n, B_{r_n^2}(x_0)\setminus B_{r_n^1}(x_1))\le\epsilon_0^2.$$ 
Applying Lemma 2.7 (iii), we can estimate that for any $\lambda>1$
\begin{eqnarray}\label{22-norm-neck31}
\big\|\nabla v_n\big\|^2_{L^{2}(B_{\lambda^{-1}}\setminus B_{\lambda r_n})}
&\le& C\Big[(A_0(r_n)+r_n)|\ln (r_n)|+C\lambda^{-2}\|\nabla v_n\|^2_{L^2(B_1\setminus B_{r_n})}\nonumber
\\
&&\quad+ \big\|\tau(v_n)\big\|^2_{M^{1,\delta}(B_1)}\Big],
\end{eqnarray}
where
$$A_0(r_n)=\inf_{h\in C^\omega(B_{2r_n})}\big\|\mathcal H(v_n)-h\|_{L^1(B_{2r_n})}
+r_n^{2-\delta}\big\|\tau(v_n)\big\|_{M^{1,\delta}(B_1)}.$$
Applying Theorem 2.8, we can estimate
\begin{eqnarray*}
\inf_{h\in C^\omega(B_{2r_n})}\big\|\mathcal H(v_n)-h\|_{L^1(B_{2r_n})}
\le C_m \Big(r_n^{2-\delta}\big\|\tau(v_n)\big\|_{M^{1,\delta}(B_1)}\Big)^{3^{1-m}},
\end{eqnarray*}
so that
\begin{equation}\label{A0-est}
A_0(r_n)\le C_m \Big(r_n^{2-\delta}\big\|\tau(v_n)\big\|_{M^{1,\delta}(B_1)}\Big)^{3^{1-m}}
\le \Big(r_n^{2-\delta}\big\|\tau(v_n)\big\|_{M^{1,\delta}(B_1)}\Big)^{3^{1-m}}
\le Cr_n^{\frac{2-\delta}{3^{m-1}}}.
\end{equation}
Substituting \eqref{A0-est} into \eqref{22-norm-neck31} and rescaling  the resulting inequality, we obtain
that for any $\lambda>1$, it holds
\begin{eqnarray}\label{22-norm-neck32}
\big\|\nabla u_n\big\|^2_{L^{2}(B_{\lambda^{-1}r_n^2}(x_1)\setminus B_{\lambda r_n^1}(x_1))}
&\le& C\Big[(r_n^{\frac{4-2\delta}{3^{m-1}}}+r_n)|\ln r_n|+
C\lambda^{-2}\|\nabla u_n\|^2_{L^2(B_{r_n^2}(x_1)\setminus B_{r_n^1}(x_1))}\nonumber\\
&&\quad+ r_n^{4-2\delta} \big\|\tau(u_n)\big\|^2_{M^{1,\delta}(\Sigma)}\big)\Big].
\end{eqnarray}
Since $\displaystyle\lim_{n\rightarrow\infty}(r_n^{\frac{2-\alpha}{3^{m+1}}}+r_n^\frac12)|\ln r_n|=0$,
it is easy to see that
\begin{eqnarray}\label{22-norm-neck33}
\lim_{n\rightarrow\infty}\big\|\nabla u_n\big\|^2_{L^{2}(B_{\lambda^{-1}r_n^2}(x_1)\setminus B_{\lambda r_n^1}(x_1))}
&\le& C\lambda^{-2},
\end{eqnarray}
and hence 
\begin{equation}\label{22-norm-vanish1}
\lim_{\lambda\rightarrow\infty}
\lim_{n\rightarrow\infty}\big\|\nabla u_n\big\|^2_{L^{2}(B_{\lambda^{-1}r_n^2}(x_1)\setminus B_{\lambda r_n^1}(x_1))}=0.
\end{equation}
(b) The neck region between the body region and a bubbling region 
$A_n=B_{\eta}(x_1)\setminus B_{r_n^1}(x_1)$ for a small $\eta>0$ and $r_n^1\rightarrow 0$.
In this case, define $r_n=\frac{r_n^1}{\eta}$ and $v_n(x)=u_n(x_1+\eta x): B_1\to N$.
Then $v_n$ is a sequence of approximate harmonic maps satisfying
$$E(v_n, B_1)\le m\epsilon_0^2;\ \ E(v_n, B_1\setminus B_{r_n})\le \epsilon_0^2.$$
Then, similar to (i), we can estimate that for any $\eta>1$, there holds
\begin{eqnarray*}
\big\|\nabla u_n\big\|^2_{L^2(B_{\lambda^{-1}\eta}(x_1)\setminus B_{\lambda r_n^1}(x_1))}
&\le& C\Big[(r_n^{\frac{4-2\delta}{3^{m-1}}}+r_n)|\ln r_n|+
C\lambda^{-2}\|\nabla u_n\|^2_{L^2(B_{\eta}(x_1)\setminus B_{r_n^1}(x_1))}\nonumber\\
&&\quad+ r_n^{4-2\delta} \big\|\tau(u_n)\big\|^2_{M^{1,\delta}(\Sigma)}\big)\Big].
\end{eqnarray*}
This again shows that 
\begin{equation}\label{22-norm-vanish2}
\lim_{\lambda\rightarrow\infty}
\lim_{n\rightarrow\infty}\big\|\nabla u_n\big\|^2_{L^{2}(B_{\lambda^{-1}\eta}(x_1)\setminus B_{\lambda r_n^1}(x_1))}=0.
\end{equation}
This establishes \eqref{bubbling2} and \eqref{bubble2.0}. Hence (ii) is proven.

To prove (iii), first observe that the uniform integrability condition \eqref{equi-cont} implies that
$$\sup_{n}\big\|\tau(u_n)\big\|_{L\log L(\Sigma)}<\infty.$$
Thus, as in the cases (a) and (b) above, we can estimate the $L^{2,1}$-norms of $\nabla u_n$ and
$L^1$-norms of $\nabla^2 u_n$ in
the neck regions by Lemma 2.7 (iv) and (v). For simplicity, let's only consider the first case, i.e., the neck region $A_n=B_{r_n^2}(x_1)\setminus B_{r_n^1}(x_1)$ given by (a) above. Then we have
\begin{eqnarray}\label{21-norm-neck31}
&&\max\Big\{\big\|\nabla v_n\big\|_{L^{2,1}(B_{\lambda^{-1}}\setminus B_{\lambda r_n})}, 
\big\|\nabla^2 v_n\big\|_{L^{1}(B_{\lambda^{-1}}\setminus B_{\lambda r_n})}\Big\}\nonumber\\
&&\le C\Big[(A_0^\frac12(r_n)+r_n^\frac12)|\ln r_n|+(C(\lambda)+\lambda^{-1})\|\nabla v_n\|_{L^2(B_1\setminus B_{r_n})}+ \big\|\tau(v_n)\|_{L\log L(B_1)}\Big].
\end{eqnarray}
This, after scaling back to $u_n$, yields
\begin{eqnarray}\label{21-norm-neck32}
&&\max\Big\{\big\|\nabla u_n\big\|_{L^{2,1}(B_{\lambda^{-1}r_n^2}(x_1)\setminus B_{\lambda r_n^1}(x_1))},
\big\|\nabla^2 u_n\big\|_{L^{1}(B_{\lambda^{-1}r_n^2}(x_1)\setminus B_{\lambda r_n^1}(x_1))}\Big\}\\
&\le& C\Big[(r_n^{\frac{2-\alpha}{3^{m+1}}}+r_n^\frac12)|\ln r_n|+
(C(\lambda)+\lambda^{-1})\|\nabla u_n\|_{L^2(B_{r_n^2}(x_1)\setminus B_{r_n^1}(x_1))}+ \big\|\tau(u_n)\|_{L\log L(B_{r_n^2}(x_1))}\Big].\nonumber
\end{eqnarray}
From \eqref{equi-cont} and \eqref{c-lambda}, we know that
$$\lim_{n\rightarrow\infty}\big\|\tau(u_n)\|_{L\log L(B_{r_n^2}(x_1))}
=0,\ \ \lim_{\lambda\rightarrow\infty} C(\lambda)=0.$$
Thus, after sending $n\rightarrow\infty$ first and then $\lambda\rightarrow\infty$ in
\eqref{21-norm-neck32}, we obtain
\begin{equation}\label{21-vanishing-neck}
\lim_{\lambda\rightarrow\infty}\lim_{n\rightarrow\infty}
\max\Big\{\big\|\nabla u_n\big\|_{L^{2,1}(B_{\lambda^{-1}r_n^2}(x_1)\setminus B_{\lambda r_n^1}(x_1))},
\big\|\nabla^2 u_n\big\|_{L^{1}(B_{\lambda^{-1}r_n^2}(x_1)\setminus B_{\lambda r_n^1})(x_1)}\Big\}=0.
\end{equation}
Note that \eqref{21-vanishing-neck} also holds in the case (b), i.e., replacing
$B_{\lambda^{-1}r_n^2}(x_1)\setminus B_{\lambda r_n^1}(x_1)$ by $B_{\lambda^{-1}\eta}(x_1)\setminus B_{\lambda r_n^1}(x_1)$ for any small $\eta>0$.

It follows from \eqref{strong-convergence} that
\begin{equation}\label{strong-convergence1}
\nabla u_n\rightarrow\nabla u \ {\rm{in}}\ L^{2,1}_{\rm{loc}}(\Sigma\setminus\mathcal S).
\end{equation}
While \eqref{bubble-convergence} implies that 
\begin{equation}\label{bubble-convergence1}
\nabla(u_n(x_n^i+r_n^i \cdot))
\rightarrow \nabla\omega_i \ {\rm{in}}\ L^{2,1}_{\rm{loc}}(\R^2), \  i=1,\cdots, m.
\end{equation}
It is clear that the $L^{2,1}$-norm identity (\ref{bubbling3}) follows from (\ref{21-vanishing-neck}),
(\ref{strong-convergence1}), and \eqref{bubble-convergence1}.
It is well-known that the oscillation convergence (\ref{oscillation-convergence}) follows from the vanishing
of $L^{2,1}$-norm of $\nabla u_n$ in the neck region and the embedding theorem (cf. \cite{tartar}):
if $f\in C^\infty_0(\R^2)$ satisfies $\nabla f\in L^{2,1}(\R^2)$, then $f\in L^\infty(\R^2)$ and
$$\|f\|_{L^\infty(\R^2)}\le C\big\|\nabla f\big\|_{L^{2,1}(\R^2)}.$$
The proof of Theorem 1.2 is now complete.
\qed

\bigskip
\noindent{\bf Proof of Corollary 1.3}: We use the same notations as in the proof of Theorem 1.2.
From (\ref{21-vanishing-neck}), to prove the $W^{2,1}$-identity (\ref{bubbling4})  it suffices to show that 
\begin{equation}\label{W21-convergence1}
\nabla^2 u_n\rightarrow \nabla^2 u \ {\rm{in}}\ L^1_{\rm{loc}}(\Sigma\setminus\mathcal S),
\end{equation}
and
\begin{equation}\label{W21-convergence2}
\nabla^2(u_n(x_n^i+r_n^i\cdot))\rightarrow \nabla^2 \omega_i \ {\rm{in}}\ L^1_{\rm{loc}}(\R^2)
\ {\rm{for}}\ i=1,\cdots, m.
\end{equation}

For \eqref{W21-convergence1}, observe that by (\ref{strong-convergence}) we have that
for any ball $B_R(x_0)\subset\subset\Sigma\setminus\mathcal S$,
$$\nabla u_n\rightarrow\nabla u\  {\rm{in}}\ L^q(B_R(x_0)), \ \forall\ 2<q<\frac{\delta}{\delta-1}.$$
From the approximate harmonic map equation \eqref{HM1}, we have
$$-\Delta(u_n-u)=(A(u_n)(\nabla u_n, \nabla u_n)-A(u)(\nabla u,\nabla u))+(\tau(u_n)-\tau),
\ {\rm{in}}\ B_R(x_0).$$
Hence, by a result of \cite{stein}, we can estimate
\begin{eqnarray*}
&&\big\|\nabla^2(u_n-u)\big\|_{L^1(B_{\frac{R}2}(x_0))}\\
&&\le C(R)\Big[\|u_n-u\|_{L^1(B_R(x_0))}+\big\|(A(u_n)(\nabla u_n, \nabla u_n)-A(u)(\nabla u,\nabla u))\big\|_{L\log L(B_R(x_0))}\\
&&\ +\big\|\tau(u_n)-\tau\big\|_{L\log L(B_R(x_0))}\Big]\\
&&\le C(R)\Big[\big\|u_n-u\big\|_{W^{1,q}(B_R(x_0))}+\big\|\tau(u_n)-\tau\big\|_{L\log L(B_R(x_0))}\Big]
\rightarrow 0,
\end{eqnarray*}
as $n\rightarrow\infty$, for some $q\in \big(2,\frac{\delta}{\delta-1}\big)$. This yields \ref{W21-convergence1})
by a simple covering argument.

\eqref{W21-convergence2} can be proven similarly. In fact, for any $R>0$, 
$v_n^i(x)=u_n(x_n^i+r_n^ix): B_R\to N$ is
an approximate harmonic map with tension field $\tau(v_n)(x)=(r_n^i)^2 \tau(u_n)(x_n^i+r_n^ix)$. Moreover,
it follows from (\ref{equi-cont}) that
$$\big\|\tau(v_n^i)\big\|_{L\log L(B_R)}\le \big\|\tau(u_n)\big\|_{L\log L(B_{Rr_n^i}(x_n^i))}\rightarrow 0,
\ {\rm{as}}\ n\rightarrow\infty.$$
Thus we can show that for any $1\le i\le m$, and for some $2<q<\frac{\delta}{\delta-1}$,
$$\big\|\nabla^2(v_n^i-\omega_i)\big\|_{L^1(B_{\frac{R}2})}
\le C(R)\Big[\big\|v_n^i-\omega_i\big\|_{W^{1,q}(B_R)}+\big\|\tau (v_n^i)\big\|_{L\log L(B_R)}\Big],$$
this, combined with (\ref{bubble-convergence}), implies (\ref{W21-convergence2}). The proof is now complete.
\qed

\bigskip
\noindent{\bf Proof of Theorem 1.1}. The proof is by a contradiction argument, 
which is similar to \cite{Lamm-Sharp} Theorem 2.1. Suppose that there  were
a sequence of approximate harmonic maps $u_n:\Sigma\to N$, with tension fields
$\tau(u_n)$, satisfying
\begin{equation}\label{global_bound}
\|\nabla u_n\|_{L^2(\Sigma)}+\|\tau(u_n)\|_{M^{1,\delta}(\Sigma)}
+\|\tau(u_n)\|_{L\log L(\Sigma)}\le \Lambda<\infty,
\end{equation}
and 
\begin{equation}\label{unbound}
\|\nabla u_n\big\|_{L^{2,1}(B_R(x_0))}+\|\nabla^2 u_n\|_{L^1(B_{R}(x_0))}\rightarrow \infty, \ {\rm{as}}\ n\rightarrow\infty,
\end{equation}
for some ball $B_R(x_0)\subset\Sigma$.

It follows from the proof of Theorem 1.2 that we can decompose $\Sigma$ into a union of body region
$\mathcal A=\Sigma\setminus\mathcal S_\eta$ (here $\mathcal S_\eta=\{x\in\Sigma:
{\rm{dist}}(x,\mathcal S)\le\eta\big\}$), bubble region $\mathcal B$,
and neck region $\mathcal N$.
From proof of Theorem 1.2, we know that \\
(i) both $\|\nabla u_n\|_{L^{2,1}(\mathcal A)}$ and $\|\nabla^2 u_n\|_{L^1(\mathcal A)}$
are uniformly bounded.\\
(ii) both $\|\nabla u_n\|_{L^{2,1}(\mathcal B)}$ and $\|\nabla^2 u_n\|_{L^1(\mathcal B)}$
are uniformly bounded.\\
Moreover, it follows from the estimate \eqref{21-norm-neck32} and the
bound \eqref{global_bound} that \\
(iii) on the neck region $\mathcal N$,
both $L^{2,1}$-norms of $\nabla u_n$ and $L^1$-norms of $\nabla^2 u_n$ are also 
uniformly bounded. 

It is clear that (i), (ii), and (iii) contradicts (\ref{unbound}). The proof is complete. \qed

\bigskip
\noindent{\bf Acknowledgement}. The author is partially supported by NSF 1522869. He wishes
to thank Zhifei Zhang for useful discussions on \cite{WWZ} and comments on this note.

{}
\end{document}